%% file: sn-article.tex
\documentclass[sn-mathphys]{sn-jnl}

\makeatletter 
\def\@acknow{}%
\long\def\EarlyAcknow#1 \par{%
\def\@acknow{\abstractfont\abstracthead*{Acknowledgments}
#1\par}}%

\def\printabstract{\ifx\@acknow\empty\else\@acknow\fi\par%
    \ifx\@abstract\empty\else\@abstract\fi\par}
\makeatother

\usepackage{algorithm,algpseudocode,amsmath, amssymb, amscd, graphicx,enumitem, wrapfig, tikz,multicol,bm, numprint, placeins,longtable,array,tabu,pgfplots,float,url,makecell,diagbox,tabularx,mathtools,hyperref,multirow,courier,subcaption,bbm,xcolor,centernot}
\usepackage{rotating}

\usepackage{pdflscape}
\usepackage{blkarray}
\usepackage{tabstackengine}
\usepackage{array}
\usepackage{booktabs}
\stackMath
\newcommand{\bla}{\color{black}}
\newcommand{\blu}{\color{black}}

\newcommand{\gre}{\color{red}}

\usepackage{amsfonts}%
\usepackage{amsthm}%
\usepackage{mathrsfs}%
\usepackage[title]{appendix}%
\usepackage{textcomp}%
\usepackage{manyfoot}%
\usepackage{listings}%

\theoremstyle{thmstyleone}%
\newtheorem{theorem}{Theorem}
\newtheorem{proposition}[theorem]{Proposition}%
\newtheorem{lem}{Lemma}

\theoremstyle{thmstyletwo}%
\newtheorem{example}{Example}%

\theoremstyle{thmstylethree}%
\newtheorem{definition}{Definition}%

\begin{document}


\title[Tight Mixed-Integer Optimization Formulations for Prescriptive Trees]{Tight Mixed-Integer Optimization Formulations for Prescriptive Trees}


\author*[1]{\fnm{Max} \sur{Biggs}}\email{biggsm@darden.virginia.edu}

\author[2]{\fnm{Georgia} \sur{{Perakis}}}
\email{georgiap@mit.edu}

\affil*[1]{\orgdiv{Darden School of Business}, \orgname{University of Virginia},\orgaddress{ \city{Charlottesville}, \state{VA}}}

\affil[2]{\orgdiv{Sloan School of Management}, \orgname{Massachusetts Institute of Technology},\orgaddress{ \city{Cambridge}, \state{MA}}}


\abstract{We focus on modeling the relationship between an input feature vector and the predicted outcome of a trained decision tree using mixed-integer optimization. This can be used in many practical applications where a decision tree or a tree ensemble is incorporated into an optimization problem to model the predicted outcomes of a decision. We propose novel tight mixed-integer optimization formulations for this problem. Existing formulations can be shown to have linear relaxations that have fractional extreme points, even for the simple case of modeling a single decision tree or a very large number of constraints, which leads to slow solve times in practice. A formulation we propose, based on a projected union of polyhedra approach, is ideal (i.e., the extreme points of the linear relaxation are integer when required) for a single decision tree. Although the formulation is generally not ideal for tree ensembles, it generally has fewer extreme points, leading to a faster time to solve. We also study formulations with a binary representation of the feature vector and present multiple approaches to tighten existing formulations. We show that fractional extreme points are removed when multiple splits are on the same feature. At an extreme, we prove that this results in ideal formulations for tree ensembles modeling a one-dimensional feature vector. Building on this result, we also show that these additional constraints result in significantly tighter linear relaxations when the feature vector is low dimensional. }




\maketitle

\section{Introduction}\label{sec1}

A fundamental problem in operations research and management science is decision-making under uncertainty. Recently, attention has been given to modeling uncertain outcomes using machine learning functions, trained from previous decisions made under a variety of circumstances \citep{bertsimas2016analytics, cheng2017maximum, tjeng2017evaluating, boob2022complexity, anderson2018strong, bunel2018unified, fischetti2018deep,kumar2019equivalent, veliborsPaper, biggs2018optimizing, bergman2022janos}. Due to the complex nature of real-world decision-making, often the model that best represents the outcomes observed is nonlinear, such as a neural network or a tree ensemble. This leads to a potentially complex optimization problem for the decision-maker to find the best decision, as predicted  by the machine learning function.
 
 An example of this occurs in reinforcement learning, where the future reward resulting from a decision is uncertain but can be approximated using machine learning models, such as decision trees or tree ensembles. In some applications, such as playing Atari video games \citep{mnih2015human}, the decision set is small so all the  decisions can be enumerated and evaluated. In comparison, in many real-world operational problems – for example, dynamic vehicle routing problems \citep{godfrey2002adaptive1, bent2007waiting, pillac2011dynamic} or kidney transplantation \citep{sonmez2017market, ashlagi2018effect}– complex decisions whose outcomes are uncertain need to be made at every stage of an online process. These decisions are often high dimensional or combinatorial in nature and subject to constraints on what is feasible. This can result in a very large action space. As a result, enumeration is no longer a tractable option, and a more disciplined optimization approach must be taken. Furthermore, the selection of the best action is further complicated by the nonlinear value function approximation.

One approach to finding optimal decisions when the outcome is estimated using a complex machine learning method is to use mixed-integer optimization (MIO) to model this relationship. In particular, there has recently been significant interest in modeling trained neural networks,  by encoding these relationships using auxiliary binary variables and constraints \citep{cheng2017maximum, tjeng2017evaluating, anderson2018strong, bunel2018unified, fischetti2018deep,kumar2019equivalent, wang2021two}. Another popular and powerful approach for supervised learning, yet one that is less studied in the prescriptive setting, is tree ensemble methods. 
\cite{veliborsPaper} provides unconstrained optimization examples in drug discovery, where a tree ensemble predicts a measure of the activity of a proposed compound, and customized price optimization, where a tree ensemble predicts the profit as a function of prices and store-level attributes. \cite{biggs2018optimizing} provide examples in real estate development of maximizing the sale price of a new house that is predicted as a function of construction decisions and location features, and a method for creating fair juries based on jurors' predicted a priori propensities to vote guilty or not due to their demographics and beliefs. These applications have nontrivial constraints, but can be represented as polyhedra with integer variables. Additional applications of trained decision trees or tree ensembles embedded in an optimization problem include retail pricing \citep{ferreira2015analytics}, assortment optimization \citep{chen2019use, chen2021assortment, chen2022decision}, last-mile delivery \citep{liu2021time}, optimal power flow \citep{halilbavsic2018data}, auction design \citep{verwer2017auction}, constraint learning \citep{maragno2021mixed} and Bayesian optimization \citep{thebelt2021entmoot}.

The goal in these works is often to propose tractable optimization formulations, which allow large problem instances to be solved in a reasonable amount of time. An important consideration when formulating these mixed-integer optimization formulations is how \textit{tight}, or strong, the formulation is. Most methods for optimizing mixed-integer formulations involve relaxing the integrality requirements on variables and solving a continuous optimization problem. In the popular branch and bound algorithm, if an optimal solution is fractional for integer variables, then multiple subproblems are created with added constraints to exclude the fractional solution. If there are fewer fractional solutions for the relaxed problem, corresponding to a tighter formulation, this can result in a significantly faster time to solve.
Furthermore, some problems can be formulated in such a way that the linear relaxation doesn't have any fractional extreme points, known as an \textit{ideal} formulation. Oftentimes, these ideal formulations can be solved extremely quickly.

Another benefit of stronger formulations is that the linear optimization (LO) relaxations provide tighter \blu dual \bla bounds, which are also useful in many applications. An example of this is evaluating the robustness of a machine learning model \citep{carlini2017towards, dvijotham2018training}. If an input can be perturbed by a practically insignificant amount and result in a significantly different prediction, this suggests that the model is not robust. Evaluating robustness can be formulated as a constrained optimization problem over local inputs to find the maximally different output. As finding the exact optimal bound can be time-consuming, an upper bound \blu on the absolute change in the objective is sufficient.\bla

\subsection{Contributions}

We model the relationship between the input feature vector and the predicted output for a trained decision tree. This can be used in a range of optimization applications involving decision trees or tree ensembles. We present a novel mixed-integer optimization formulation based on a projected \textit{union of polyhedra} approach, which we prove is ideal for a single tree and has fewer constraints and variables than existing formulations. 
We show existing mixed-integer optimization formulations for modeling trees either are not ideal for a single tree \citep{biggs2018optimizing, veliborsPaper}; or contain significantly more constraints \citep{kim2022reciprocity}\footnote{Our results were developed independently from this recent paper (and exist in an earlier version of this paper from 2020, \cite{biggs2020dynamic})}, leading to substantially slower times to solve in practice. 
Our formulation applies to general feature vectors compared to \cite{veliborsPaper, kim2022reciprocity}, which use binary encodings of the feature vector and are more difficult to incorporate into a constrained optimization formulation.  
While the formulation we present is generally not ideal when we impose polyhedral constraints on the decision, or when multiple trees are used in an ensemble model, the formulation generally excludes fractional extreme points present in \cite{biggs2018optimizing} and \cite{veliborsPaper}, leading to tighter formulations. 

We also present new formulations that use a binary feature vector representation, as proposed in \cite{veliborsPaper}. Despite the aforementioned difficulties with constrained optimization formulations, these formulations do appear to have some advantages regarding the branching behavior in the MIO solver, leading to a faster time to solve in some instances. We propose different constraints that can be added to tighten the formulation from \cite{veliborsPaper}. The \textit{expset} formulation is based on exploiting the greater than or equal to representation of the feature vector from \cite{veliborsPaper}, leading to larger groups of leaf variables being turned off when a split is made. The \textit{elbow} formulation removes specific fractional solutions that arise when there are nested branches on the same feature in a tree. We characterize the conditions in which each of these constraints removes fractional solutions, which generally occurs in scenarios where there are multiple splits on the same feature. Extending this, we show that the \textit{expset} formulation leads to an ideal formulation when all the splits are on the same feature, which occurs for tree ensembles when the feature vector is one-dimensional. This property doesn't hold for the formulations in \cite{veliborsPaper,chen2021assortment,kim2022reciprocity}, and can be contrasted with the results in \cite{chen2021assortment, kim2022reciprocity}, which present ideal formulations for a single decision tree (but with many dimensions) for a binary encoded feature vector. 
These results provide insights for the practitioner on when different formulations might be tighter. When there are many trees in the ensemble but relatively few variables, the \textit{expset} formulation is likely to be tighter. When there are few trees but many variables, the \textit{union of polyhedra} or formulation from \cite{kim2022reciprocity}  is likely to be tighter.

We explore the performance of these approaches through extensive simulations. In agreement with our theoretical findings, we show that in many instances, the \textit{union of polyhedra} formulation appears to have significant solve time improvements for tree ensembles with few but large trees. Similarly, the \textit{elbow} offers improvements for problems with few features. Despite the theoretical appeal of the tightness of \cite{kim2022reciprocity}, we show that in practice, it is much slower than the other proposed approaches due to the very large number of constraints added. While the \textit{expset} formulation generally doesn't offer faster solve times, we show that the linear relaxations it provides can be significantly stronger. This is useful in many applications where a bound on the optimal solution is desired, particularly for trees with few features.

\section{Preliminaries}


Given a feature vector $\bm{w} \in \mathcal{D} \subseteq \mathbb{R}^d$, our goal is to model the output $y_t \in \mathcal{Y} \subseteq \mathbb{R}$ of a decision tree $f^{(t)}(\bm{w})$ using a mixed-integer optimization formulation, where $t$ corresponds to the index of the tree in the case of a tree ensemble. More formally, we model the graph, $gr(f^{(t)};\mathcal{D}) = \{\bm{w},y_t |\bm{w} \in \mathcal{D}, y_t=f^{(t)}(\bm{w})\}$. With such a formulation, we can easily model a range of practical applications, such as finding an optimal feature vector to maximize the predicted outcome of a tree ensemble, $\max_{w \in \mathcal{D}} \sum_{t=1}^T f^{(t)}(\bm{w})$, or solving a reinforcement learning subproblem with complex constraints where the value function is given by a decision tree. 


\subsection{Decision trees}

\begin{figure}[]
  \centering
  \begin{subfigure}{0.48\textwidth}
    \centering
	\includegraphics[width=0.7\textwidth]{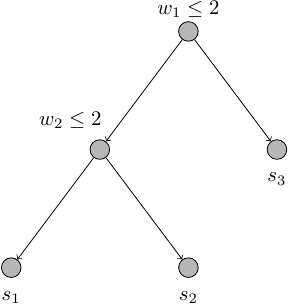}
	\centering
	\caption{Example decision tree}
	\label{fig:union_tree}
  \end{subfigure}
    \begin{subfigure}{0.48\textwidth}
    \centering
	\includegraphics[width=0.8\textwidth]{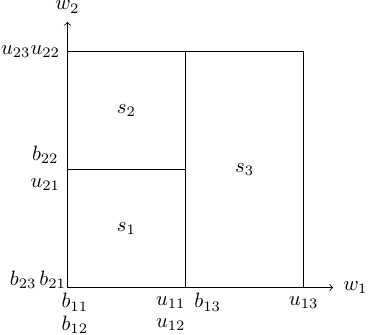}
	\centering
	\caption{Partition of feature space with bounds}
	\label{fig:union_partition}
  \end{subfigure}
  \caption{Examples of decision tree with corresponding notation and partition of the feature space}
  \label{fig:simple_dec_tree}
\end{figure}


A decision tree $f^{(t)}(\bm{w})$ with $p$ leaves is a piecewise constant function, where a constant outcome $s_l \in \mathcal{Y} \subseteq \mathbb{R}$ is predicted if feature vector $\bm{w}$ falls within a particular leaf $\mathcal{L}_l, l \in [p]$, so that $f^{(t)}(\bm{w})=s_l ~\text{if}~ \bm{w} \in  \mathcal{L}_l$. Each leaf, $\mathcal{L}_l$, is a hyperrectangular set defined by an upper $u_{li}$ and a lower (bottom) $b_{li}$ bound for each feature dimension $w_i, i  \in~ [d]$. Throughout, we assume $\mathcal{D}$, and therefore $w_i$, is bounded. 

\vspace{10pt}
\begin{definition} A leaf in a decision tree satisfies:
\label{leaf_def}
\begin{subequations}
\begin{align}
\mathcal{L}_l=\{ \bm{w} , y~  |~ w_i &\leq u_{li} \qquad ~\forall ~i  \in~ [d],\\
w_i &\geq b_{li} \qquad ~\forall ~i  \in~ [d],\\
y&=s_l \}. 
\end{align}
\end{subequations}
\end{definition}

A hierarchy of axis-aligned splits defines the upper bounds and lower bounds associated with each leaf, each of which is on a single variable, i.e., $w_i \leq \theta$. These splits define the tree and partition the feature space into leaves. We denote $\textbf{splits}(t)$ as the set of splits corresponding to tree $t \in T$, $\textbf{left}(s)$ as the set of leaves to the left of split $s$ in the tree (i.e., those that satisfy the split condition  $w_i \leq \theta$), and $\textbf{right}(s)$ as the set of leaves to the right for which  $w_i \geq \theta$. The upper bounds $u_{li}$ are defined by the threshold of the left splits that lead to the leaf, while the lower bounds $b_{li}$  are defined by the thresholds of the right splits.\footnote{ We note that our definition of a leaf differs slightly from the standard definition of a leaf used in a decision tree, where there is typically a strict inequality associated with a threshold (i.e., the lower bound leaf would be defined by $w_i > \theta$). We use our definition ($w_i \geq \theta$) due to the inability of mixed integer optimization to model open sets. As such, if there is a vector precisely at the threshold $w_i = \theta$, it could be in either leaf, but when maximized/minimized in an optimization context, $\bm{w}$ will end up being in the leaf with the higher/lower predicted outcome.} In the case where there are multiple axis-aligned splits along a dimension leading to a leaf (i.e., $w_1 \leq 5$ then $w_1 \leq 2$), the upper bound will be the minimum of all less than splits, while the lower bound will be the maximum. When there are no splits on a feature, the upper and lower bounds on the leaf are the upper and lower bounds on the feature vector.

\subsection{Mixed-integer optimization}
We aim to model the graph $gr(f; \mathcal{D})$ using mixed-integer optimization. To facilitate this, often auxiliary continuous $\bm{q} \in \mathbb{R}^n$ and integer variables are introduced to help model the complex relationships between variables, although the formulations we study require only binary variables $\bm{z} \in  \{0,1\}^m$. A mixed-integer optimization formulation consists of linear constraints on $(\bm{w},y,\bm{q},\bm{z}) \in \mathbb{R}^{d+1+n+m}$ which define a polyhedron $Q$, combined with binary constraints on $\bm{z} \in \{0,1\}^m$. For a valid formulation, the set $(\bm{w},y)$ associated with a feasible solution $(\bm{w},y,\bm{q},\bm{z}) \in Q \cap \mathbb{R}^{d+1+n} \times \{0,1\}^m $ must be the same as the graph we desire to model $(\bm{w},y) \in gr(f; \mathcal{D})$. More formally, the auxiliary variables ($\bm{q},\bm{z}$) are removed via an orthogonal projection $Proj_{\bm{w},y}(Q)=\{\bm{w},y~|~ \exists~ \bm{q}, \bm{z} ~s.t.~ \bm{w},y,\bm{q},\bm{z} \in Q \}$, to leave a set of feasible $(\bm{w},y)$. Therefore, a valid mixed-integer optimization formulation may be defined as:
\vspace{10pt}
\begin{definition} A valid mixed-integer optimization formulation satisfies:
$$ gr(f; \mathcal{D})= Proj_{\bm{w},y} (Q \cap \mathbb{R}^{d+1+n} \times \{0,1\}^m).$$
\end{definition}

We will refer to $Q$ as the linear relaxation of the formulation, which is the MIO formulation with the integrality requirements removed. A MIO formulation is ideal if the extreme points of the polyhedron are binary for those variables that are required to be: 
\vspace{10pt}
\begin{definition} An ideal formulation satisfies: $$\text{ext}(Q) \subseteq \mathbb{R}^{d+1+n} \times \{0,1\}^m$$ 
where $\text{ext}(Q)$ are the extreme points of the polyhedron $Q$.
\end{definition}








\section{Further relevant literature}


 As previously mentioned, modeling trained tree ensembles using mixed-integer optimization is studied in  \cite{biggs2018optimizing,veliborsPaper, kim2022reciprocity}. \cite{veliborsPaper} proved this problem is NP-Hard and proposed formulations for unconstrained optimization problems or problems with simple box constraints on each variable. \cite{mistry2021mixed} provide a customized branch and bound algorithm for optimizing gradient-boosted tree ensembles based on the MIO formulation in \cite{veliborsPaper}, while \cite{perakis2021motem} also propose a customized branching procedure. \cite{biggs2018optimizing} proposes formulations that include polyhedral constraints. This approach uses the big-M approach to linearize the nonlinear behavior of the trees. To optimize large tree ensembles in a reasonable amount of time, both \cite{veliborsPaper} and \cite{biggs2018optimizing} offer ways to decompose a large tree ensemble and propose heuristic approaches that involve truncating trees to a limited depth \citep{veliborsPaper} or sampling a subset of the trees \citep{biggs2018optimizing}. \cite{kim2022reciprocity} highlights equivalences between tree ensemble optimization and multilinear optimization and provides formulations based on techniques from multilinear optimization. As previously mentioned, these formulations are ideal for a single tree but introduce many constraints, so the time to solve is often significantly longer. The formulations in \cite{kim2022reciprocity} generalize those in \cite{chen2021assortment}. \cite{chen2021assortment} studies the assortment optimization setting, where each feature is binary and therefore branched on at most once in each decision tree (corresponding to the inclusion of a product in an assortment or not). All these approaches involve solving a mixed-integer optimization formulation of an ensemble of trees. 

 \blu
 There also exists a rich literature on the related but distinct problem of training decision trees using mixed-integer optimization (see, for example, \cite{bertsimas2017optimal,michini2024polyhedral,aghaei2024strong}), rather than our problem of finding the optimal decision, given an already trained decision tree, or ensemble. \bla

  

\subsection{Formulation from \cite{veliborsPaper}}
We review the formulation from \cite{veliborsPaper} both as a benchmark and to motivate the formulations we propose. Rather than linking the feature vector $\bm{w}$ directly to the output $f^t(\bm{w})$, \cite{veliborsPaper} uses a binary representation of the feature vector $\bm{w}$, which represents whether the feature falls below each split in the tree. Specifically, binary variables are introduced with
$$
x_{ij}=
\begin{cases}
1 ~~~ \text{if} ~w_i \leq \theta_{ij}\\
0 ~~~ \text{if} ~w_i \geq \theta_{ij}
\end{cases}
$$
 where $\theta_{ij}$ is the $j^{th}$ largest split threshold associated with dimension $i$. As a result, the $\bm{x}_i$ vector has the structure of consecutive 0's, followed by consecutive 1's. For example, $\bm{x}_i=\{0,1,1\}$ would correspond to a solution that falls between the first and second thresholds. A drawback of this approach is that typically, additional constraints and variables are needed to place constraints on the input vector. 
 
 
 To introduce the formulation from \cite{veliborsPaper}, we need to introduce some additional notation. $C(s)$ corresponds to the ranking of threshold $s$ relative to the size of other thresholds for that feature, and $V(s)$ corresponds to the feature involved in the split. For example, if $\theta_{ij}$ is the $j^{th}$ largest threshold for feature $i$ associated with split $s$, then $C(s)=j$ and $V(s)=i$. $K_i$ denotes the number of thresholds for feature $i$. Auxiliary variables $\bm{z}$ are introduced, where $z_l=1$ if the feature vector falls in leaf $\mathcal{L}_l$. The polyhedron $Q^{misic}$, which links the binary representation $\bm{x}$ to the predicted outcome $y$, is:
  \begin{subequations}
  \label{velibors_formulation}
\begin{align}
     Q^{misic}  = \{ \bm{x},y,\bm{z} ~ |~ &  \sum_{l \in \textbf{left}(s)} z_{l} \leq x_{V(s)C(s)}  \qquad \forall s ~\in  ~ \textbf{splits}(t)  \label{left_constraint}\\
    & \sum_{l \in \textbf{right}(s)} z_{l} \leq 1- x_{V(s)C(s)}  \qquad \forall s ~\in  ~ \textbf{splits}(t) \label{right_constraint}\\
    &  x_{ij} \leq x_{ij+1}  \qquad \forall i ~\in  ~ [d], ~ \forall j ~\in  ~ [K_i] \label{bigger_than_cosntraint} \\
    & \sum_{l =1}^p z_{l} = 1,~~~ y = \sum_{l = 1}^p s_l z_l  \\
    & \bm{x} \in [0,1]^{K_i} \qquad \forall i \in [d], ~\bm{z} \geq 0 \}.
\end{align}
  \end{subequations}
  
 The corresponding MIO formulation imposes binary constraints on  $\bm{x}\in \{0,1\}^{K_i} ~ \forall i \in [d]$, but they are not necessary for $\bm{z}$. Constraint (\ref{left_constraint}) enforces that if the condition at a split is not satisfied, $x_{V(s)C(s)}=0$, then the solution does not fall within a leaf to the left of that split in the tree, so $z_l= 0 ~ \forall  l ~\in \textbf{left}(s)$. Conversely in constraint (\ref{right_constraint}), if the split is satisfied, $x_{V(s)C(s)}=1$, then all leaves to the right are set to 0. Constraint (\ref{bigger_than_cosntraint}) links the solution to the feature vector across trees. If the solution is less than the $j^{th}$ split, $x_{ij}=1$, then the solution must also be less than all splits greater than this. As such, $x_{ik}=1 ~ \forall j< k < K_i$, and the vector has the structure of consecutive zeros followed by consecutive ones.
  
  An issue with the formulations presented in both \cite{veliborsPaper} and \cite{biggs2018optimizing} is that the linear relaxation can have many fractional solutions. This can make the MIO slow to solve. In fact, neither formulation is ideal even for the simple case of modeling a single decision tree without any additional constraints on a feasible decision, as we show in the following example.

\begin{figure}[]
  \centering
  \begin{subfigure}{0.48\textwidth}
    \centering
	\includegraphics[width=0.7\textwidth]{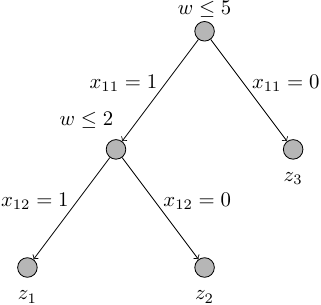}
	\centering
	\caption{\cite{veliborsPaper}}
	\label{fig:velibor_not_ideal}
  \end{subfigure}
    \begin{subfigure}{0.48\textwidth}
    \centering
	\includegraphics[width=0.7\textwidth]{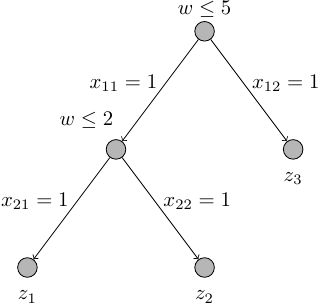}
	\centering
	\caption{\cite{biggs2018optimizing}}
	\label{fig:biggs_not_ideal}
  \end{subfigure}
  \caption{Examples of trees with fractional solutions and notation}
  \label{fig:trees_not_ideaal_current}
\end{figure}

\vspace{10pt}
\begin{example} [\cite{veliborsPaper} not ideal for a single tree with a single feature]
\label{velibor_not_ideal}
Suppose there is a tree that first branches on the condition $w \leq 5$ and then on $w \leq 2$, as shown in Figure \ref{fig:velibor_not_ideal}. In this example, $x_{11}=1$ if $w \leq 5$, and 0 otherwise, while $x_{12}=1$ if $w \leq 2$. The variables $z_l=1$ if the solution is in leaf $\mathcal{L}_l$. The resulting linear relaxation from \cite{veliborsPaper} is:
\begin{align*}
    \{ \bm{x},\bm{z}~ |~ z_2 &\leq 1-x_{12}, &  &z_3 \leq 1-x_{11}, & &x_{12}  \leq x_{11} && 0 \leq  \bm{x} \leq 1,\\
    ~ z_1  &\leq x_{12}, & & z_1 + z_2 \leq x_{11}, & & z_1+z_2+z_3 =1,
     && 0 \leq  \bm{z}  \}.
\end{align*}

This has an extreme point at $z_1=0,~ z_2=0.5,~ z_3=0.5,~ x_{11}=0.5, ~x_{12}=0.5$, when constraints $ z_2 \leq 1-x_{12},~z_3  \leq 1-x_{11}, ~x_{12}  \leq x_{11}, ~z_1+z_2+z_3 =1,~z_1 \geq 0 $ are active.  
\end{example}

\vspace{10pt}
\begin{example} [\cite{biggs2018optimizing} not ideal for a single tree with a single feature]
\label{biggs_not_ideal}
Again, suppose there is a tree that first branches on the condition $w \leq 5$ and then on $w \leq 2$, as shown in Figure \ref{fig:biggs_not_ideal}. This formulation uses a slightly different notation, where $x_{ij}=1$ if the arc is on the path to the active leaf,  $i$ corresponds to the parent node, $j=1$ refers to the left branch, and $j=2$ refers to the right branch. For example, if $w \leq 2$, then $x_{11},x_{21}=1$, while $x_{12},x_{22}=0$.  We also assume $w$ is bounded, $0 \leq w \leq 10$, and following guidance in \cite{biggs2018optimizing} for choosing the big-M value, we set $M=15$. The resulting formulation in \cite{biggs2018optimizing} is:
\begin{align*}
    \{ \bm{x},w ~ |~ &w - 15(1-x_{11}) \leq 5 , 
    && w - 15(1-x_{21}) \leq 2,
    && x_{21}+x_{22} = x_{11}, \\
    & w + 15(1-x_{12}) \geq 5 ,
    && w + 15(1-x_{22}) \geq 2 ,
    && x_{12} + x_{21}+x_{22}  = 1,\\
    & 0 \leq w \leq 10, && 0 \leq \bm{x} \leq 1 \}.
\end{align*}

This has an extreme point at $ x_{11} = 1/3, ~ x_{12} = 2/3,~  x_{21} = 1/3,~  x_{22} = 0, w=0 $, when constraints $  w + 15(1-x_{12}) \geq 5 ,~x_{21}+x_{22} = x_{11}, ~x_{11} + x_{12} + x_{21}+x_{22}  = 1, ~w\geq 0, ~x_{22}  \geq 0 $ are active. Furthermore, this is not just a consequence of the choice of $M$ but is still an issue regardless of this choice. 
\end{example}

 \section{Union of polyhedra formulation}
 
We propose an alternative MIO formulation for decision trees, which is tighter in the sense that it is ideal for modeling a single tree, unlike those presented in Example \ref{velibor_not_ideal} and \ref{biggs_not_ideal}. In contrast with the formulations in \cite{veliborsPaper}, \cite{chen2021assortment} and \cite{kim2022reciprocity}, our proposed formulation directly relates the feature vector $\bm{w}$, to the output $f^{(t)}(\bm{w})$, instead of using a binary representation of the feature vector. This has the advantage that constraints can be placed directly on the feature vector $\bm{w}$ for problems with additional constraints that need to be modeled.

\blu 
To develop our formulation, we explicitly consider the decision tree as a union of polyhedra \citep{balas1985disjunctive} corresponding to the leaf sets from Definition \ref{leaf_def}, $\cup_{l \in [p]} \mathcal{L}_l$. The leaf sets $\mathcal{L}_l$ are hyperrectangles that partition the feature space.  Leveraging established results on disjunctive formulations \citep{balas1985disjunctive}, this union can be modeled using the classical extended formulation approach from \cite{jeroslow1987representability}. This formulation, also recognized as a ``multiple-choice" formulation \citep{vielma2011modeling} or ``convex hull refromulation" \citep{grossmann2002review}, introduces auxiliary variables to explicitly capture which leaf set the solution resides in: \bla
 \begin{subequations}
 \label{union_polyhedra_formulation}
\begin{align}
Q^{ext} = \{  \bm{w}, y , \bar{\bm{w}}, \bar{\bm{y}}, \bm{z} | ~ & u_{li} z_l \geq \bar{w}_{li} \qquad \forall i \in [d], ~\forall l \in [p] \label{union_upper_bound}\\
&  b_{li} z_l \leq \bar{w}_{li}  \qquad \forall i \in [d],  ~\forall l \in [p] \label{union_lower_bound}\\
& \bar{y}_l = s_l z_l,\qquad \forall l \in [p] \label{union_score} \\
& \sum_{l = 1}^p z_l =1,  \label{union_choose_one_set}\\
& w_i = \sum_{l = 1}^p \bar{w}_{li}  \qquad \forall i \in [d] \label{union_form_w} \\ 
& y = \sum_{l = 1}^p \bar{y}_l  \label{union_form_y} \\
& z_l \in [0,1]  \qquad \forall l \in [p] \label{union_z} \}.
\end{align}
 \end{subequations}

 The formulation works by creating auxiliary copies of each variable, $\bar{\bm{w}}_l \in \mathbb{R}^d, \bar{y}_l \in \mathbb{R}$, corresponding to each leaf $l \in [p]$. With a slight abuse of notation, $\bar{\bm{w}}$ corresponds to the matrix $[\bar{\bm{w}}_1, ..., \bar{\bm{w}}_p]$. Auxiliary binary variables $\bm{z} \in \{0,1\}^p$ are also introduced, which indicate which leaf the solution falls into. When $z_l=1$, constraints (\ref{union_upper_bound}),  (\ref{union_lower_bound}), and (\ref{union_score}) define the feasible region and score for that leaf. When $z_l=0$, these constraints enforce that $\bar{\bm{w}}_l$ is set to be a vector of zeros. Constraints (\ref{union_choose_one_set}) ensure only one leaf is chosen. Constraint (\ref{union_form_w}) and (\ref{union_form_y}) in turn define $\bm{w}$ and $y$ according to which leaf is active. 



This formulation is ideal, as proved in \cite{jeroslow1984modelling} and \cite{balas1985disjunctive}, so the linear relaxation is guaranteed to have integer extreme points. However, these formulations often have computational issues when solved in practice \citep{vielma2019small}. This formulation introduces a large number of auxiliary variables ($(p+1)(d+2)$ variables in total), as well as many constraints $(2pd +3p+d+1)$. It is well known that these formulations suffer from degeneracy, as many of the auxiliary variables are set to be 0, often resulting in poor performance in practice \citep{vielma2019small}. 

\blu
This formulation can be relaxed through the aggregation of the leaf-specific constraints. Specifically, summing constraint (\ref{union_upper_bound}) across all leaves and substituting \(w_i = \sum_{l=1}^p \bar{w}_{li}\) from (\ref{union_form_w}), we  obtain:
$$\sum_{l=1}^p u_{li} z_l \geq \sum_{l=1}^p \bar{w}_{li} = w_i, \quad \forall i \in [d]. $$
Similarly, summing constraint (\ref{union_lower_bound}) yields:
$$ \sum_{l=1}^p b_{li} z_l \leq \sum_{l=1}^p \bar{w}_{li} = w_i, \quad \forall i \in [d]. $$
Finally, summing constraint (\ref{union_score}) over all leaves and substituting \(y = \sum_{l=1}^p \bar{y}_l\) from (\ref{union_form_y}) gives:
$$y = \sum_{l=1}^p s_l z_l.$$

With these aggregations, we arrive at a formulation that is a relaxation of (\ref{union_polyhedra_formulation}). This significantly smaller formulation involves only the original variables and the binary selection variables:
\bla
 
 
 
 \begin{subequations}
 \label{projected_formulation}
 \begin{align}
Q^{proj} = \{ \bm{w}, y, \bm{z} | ~ &\sum_{l = 1}^p u_{li} z_l \geq w_i \qquad \forall i \in [d], \label{upper_bound}\\
& \sum_{l = 1}^p b_{li} z_l \leq w_i  \qquad \forall i \in [d] , \label{lower_bound}\\
& y = \sum_{l = 1}^p s_l z_l  \\
& \sum_{l = 1}^p z_l =1   \label{choose_one_set}\\
& z_l \in [0,1]  \qquad \forall l \in [p]  \label{aux_binary} \}.
\end{align}
\end{subequations}

\blu 
Surprisingly, we can show that this formulation is also the projection, via Fourier-Motzkin elimination, of $Q^{ext}$ onto $\bm{w},y,$ and $\bm{z}$. Since formulation (\ref{union_polyhedra_formulation}) is ideal, the projection is also ideal. As a result, we can prove this formulation is ideal for a single tree: \bla
\vspace{10pt}
\begin{theorem}[Ideal formulation for a tree]
\label{ideal_tree}
$Proj_{\bm{w},y, \bm{z}} (Q^{ext})= Q^{\text{proj}}$. Furthermore, the polyhedron $Q^{\text{proj}}$ is ideal.
\end{theorem}
\vspace{10pt}
This is formally proved in Appendix \ref{proof_ideal}. These ideal projected formulations always exist, but in general, the projection is not a tractable operation and can result in a formulation with exponentially many constraints. In this special case, the resulting formulation (\ref{projected_formulation}) has only $2d+2$ constraints (in addition to binary constraints) and $p+d+1$ variables. Compared to formulation (\ref{union_polyhedra_formulation}), this has significantly fewer variables and therefore does not suffer from degeneracy to the same extent. 

 We also note that this formulation has considerably fewer constraints than in \cite{veliborsPaper}, which has approximately $3(\sum_{i=1}^d K_i)$ constraints (3 constraints for each split) and $\sum_{i=1}^d K_i + p$ variables. We recall that $K_i$ is the number of splits for feature $i$. For an axis-aligned tree, the total number of leaves equals the total number of splits plus one, so this corresponds to approximately $3p$ constraints and $2p$ variables. In most applications, $d <<p$, so this is substantially more than in formulation (\ref{projected_formulation}). The number of constraints is also substantially less than in \cite{kim2022reciprocity}, which is of the order $\sum_{i=1}^d K_i^2$ with $\sum_{i=1}^d K_i + p$ variables. For tree ensembles with $T$ trees, this scales as $T(\sum_{i=1}^d K_i^2)$, while in \cite{veliborsPaper} the scaling is still $3(\sum_{i=1}^d K_i)$ (although $K_i$ increases with $T$). Although this quadratic dependence may seem mild, considering MIP solve times are already exponential in the problem size, it leads to dramatically longer solve times in practice. The practical formulation sizes for various formulations are empirically studied in Table \ref{table:problem_size_small}. 

The significance of Theorem \ref{ideal_tree} is that it suggests that tree-based optimization approaches that use formulation (\ref{projected_formulation}) will be tighter than those used in \cite{biggs2018optimizing} or \cite{veliborsPaper}. Specifically, there are fractional solutions for each tree, as shown in Examples \ref{velibor_not_ideal} and \ref{biggs_not_ideal}, which do not exist in formulation (\ref{projected_formulation}). However, in general, the intersection of different tree polytopes, as occurs in tree ensemble optimization, introduces additional fractional solutions.  This also occurs for the intersection of a tree polytope and additional polyhedral constraints. However, in practice, this formulation often results in a faster time to solve, particularly for forests with relatively few trees.

If formulation (\ref{projected_formulation}) is reformulated slightly, we can prove some additional favorable properties, including, in particular, that the constraints are facet-defining for the polyhedron.
\vspace{10pt}
\begin{definition}
A face $\mathcal{F}$ of a polyhedron $\mathcal{P}$, represented by the inequality $\bm{a}'\bm{x} \geq b$, is called a facet of $\mathcal{P}$ if $dim(\mathcal{F} ) = dim(\mathcal{P} ) - 1$.
\end{definition}
\vspace{10pt}
One of the variables $z_p$ can be eliminated through the substitution $z_p=1- \sum_{l = 1}^{p-1} z_l $. Consequently, $\bm{z} \in \{0,1\}^{p-1}$ and as a result, $\bm{z}=0$ implies $\bm{w} \in \mathcal{L}_p$ (see definition \ref{leaf_def}). This leads to the following formulation:
 \begin{subequations}
  \label{facet_defining_formulation}
\begin{align}
Q^{facet} = \{  \bm{w},y,\bm{z} | ~ &u_{pi}+ \sum_{l=1}^{p-1} (u_{li}-u_{pi})z_l \geq w_i \qquad \forall i \in [d], \label{upper_bound_facet}\\
&  b_{pi} + \sum_{l=1}^{p-1} (b_{li} - b_{pi})z_l \leq w_i  \qquad \forall i \in [d] , \label{lower_bound_facet}\\
& y =  s_p + \sum_{l=1}^{p-1} z_l (s_l-s_p) \\
& \sum_{l = 1}^{p-1} z_l =1  \\
& z_l \in [0,1]  \qquad \forall l \in [p-1]  \}.
\end{align}
 \end{subequations}

We can show that under mild assumptions, (\ref{upper_bound_facet}) and (\ref{lower_bound_facet}) are facet-defining for $Q^{facet}$. 

\vspace{10pt}
\begin{lem}
\label{facet_def}
For all $l \in [p]$, assume $\mathcal{L}_l$ is non-empty and $\mathcal{L}_l$ is full dimensional, i.e., $dim(\mathcal{L}_l)= d$. Then constraints (\ref{upper_bound_facet}) and (\ref{lower_bound_facet}) are facet-defining.
\end{lem}
\vspace{10pt}
 
 This is proved in Appendix \ref{proof_facet_defining} with a proof technique similar to that in \cite{anderson2018strong}. This result is significant because it suggests there is no redundancy in formulation (\ref{facet_defining_formulation}). MIO formulations generally take longer to solve when there are redundant variables and constraints.

\subsection{Extensions to tree ensembles and additional constraints}
\label{sec:tree_ensebles}

 The formulation can be applied to tree ensembles such as random forests or gradient-boosted tree ensembles. While the polyhedron modeling an individual tree is ideal, this formulation is not ideal in general as shown in this section. An alternative, but weaker, notion of tightness is whether a formulation is sharp. For a sharp formulation, the projection of the polyhedron $Q$ onto the original variables $\bm{w},y$ is equal to the convex hull $(\text{conv}(\cdot))$ of the graph $gr(f; \mathcal{D})$. This is formalized as follows:  
 
\vspace{10pt}
 \begin{definition} A sharp formulation satisfies:
      $$\text{conv}(gr(f; \mathcal{D}))=Proj_{\bm{w},y}(Q).$$
 \end{definition}
 
An ideal formulation is also sharp, but a sharp formulation isn't necessarily ideal. In Example \ref{example_tree_not_ideal} we give a simple tree ensemble that illustrates that the \textit{union of polyhedra} formulation is not ideal and not sharp for an ensemble.


\begin{figure}[]
  \centering
  \begin{subfigure}{0.45\textwidth}
    \centering
	\includegraphics[width=0.8\textwidth]{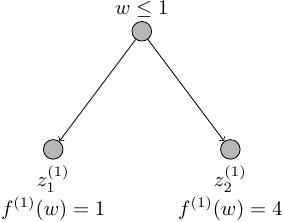}
	\caption{ Tree 1}
  \end{subfigure}
    \begin{subfigure}{0.45\textwidth}
    \centering
	\includegraphics[width=0.8\textwidth]{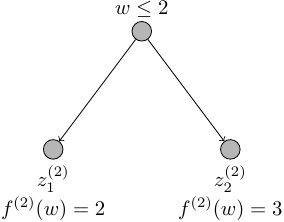}
	\caption{ Tree 2}
  \end{subfigure} \\
    \begin{subfigure}{0.8\textwidth}
    \centering
	\includegraphics[width=0.8\textwidth]{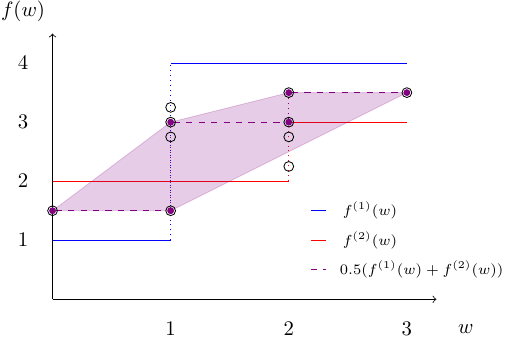}
	\caption{Graph of tree ensemble}
 \label{fig:graph_ensemble}
  \end{subfigure}
  \caption{Tree ensemble formulation is not ideal or sharp. Extreme points of $Q^{proj}$ are shown with hollow circles, while the convex hull of the tree ensemble graph is shown in shaded purple.}
  \label{example_1_fig}
\end{figure}

\vspace{10pt}
\begin{example}[Intersection of trees, using $Q^{proj}_t$ for each tree $t$, is not ideal or sharp]
\label{example_tree_not_ideal}

Suppose we have the following two trees in an ensemble:
$$f^{(1)}(w) = \begin{cases} 1 ~~~ 0 \leq w \leq 1 \\ 4 ~~~ 1 < w \leq 3 \end{cases} \qquad f^{(2)}(w)= \begin{cases} 2 ~~~ 0 \leq w \leq 2 \\ 3 ~~~ 2 < w \leq 3. \end{cases}  $$

This leads to a tree ensemble:
$$0.5(f^{(1)}(w)+f^{(2)}(w)) =\begin{cases} 1.5 ~~~ 0 \leq w \leq 1 \\ 3 ~~~~~ 1 < w \leq 2 \\ 3.5 ~~~ 2 < w \leq 3. \end{cases}$$

This is visualized in Figure \ref{example_tree_not_ideal}, where $f^{(1)}(w)$ is the blue line, $f^{(2)}(w)$ is the red line and the ensemble $0.5(f^{(1)}(w)+f^{(2)}(w))$ is the purple dashed line. The \textit{union of polyhedra} formulation for this is as follows:
  \begin{subequations}
\begin{align*}
\{ w, y, \bm{z}~|~ & z_2^{(1)} \leq w, 
~~~ z_1^{(1)} + 3z_2^{(1)} \geq w, 
  ~~~~ z_1^{(1)} + z_2^{(1)} =1 , \\
  &2z_2^{(2)} \leq w,
  ~~ 2z_1^{(2)} + 3z_2^{(2)} \geq w,
 ~~ z_1^{(2)} + z_2^{(2)} =1 , \\
  &y = 0.5\left( z_1^{(1)}+4z_2^{(1)}+2z_1^{(2)} +3z_2^{(2)} \right),
 ~~ \bm{z},\bm{w} \geq 0 \}.
\end{align*}
 \end{subequations}
A basic feasible solution for this formulation is  $w=1,~ z_1^{(1)}=0, ~z_1^{(2)} =1, ~ z_2^{(1)} =0.5,~ z_2^{(2)} =0.5,~ y=3.25$, which is not integral, so the formulation is not ideal. Furthermore, the projected solution, $w=1,~ y=3.25$, is not in the convex hull of $0.5(f^{(1)}(w)+f^{(2)}(w))$, so the formulation is not sharp. This can be observed in Figure \ref{fig:graph_ensemble}, where the convex hull of the graph of the tree ensemble is shown in shaded purple. The extreme points of  $Q^{proj}$ projected into $w,y$ space are shown with hollow circles. As can be observed, there are two extreme points of $Q^{proj}$ that lie outside the convex hull of the graph.  
 \end{example}
\vspace{10pt}
We also provide an example illustrating that adding additional constraints to the feature vector, which may be useful for many practical applications, is not ideal.
\vspace{10pt}
\begin{example}[Adding additional constraints to a tree is not ideal]
\label{const_example_tree_not_ideal}

Take the tree from Figure \ref{fig:simple_dec_tree}. Suppose that we add a simple constraint that $w_1+w_2 \leq 3$. Suppose additionally that there are upper and lower bounds on each feature, such that $0 \leq w_1, w_2\leq3$. The \textit{union of polyhedra} formulation is:
   \begin{subequations}
\begin{align*}
\{ w_1, w_2, \bm{z}~|~ && 2(z_1+z_2)+3z_3 \geq w_1,  
 && 2z_1+3(z_2+z_3) \geq w_2,  
 && z_1+z_2+z_3 =1 \\
&&2z_3 \leq w_1,  
&&2z_2 \leq w_2 ,
  && w_1+w_2 \leq 3,~  \bm{z} \geq 0 \}.
\end{align*}
 \end{subequations}

 This has a fractional solution $w_1 = 2/3,~ w_2= 7/3,~z_1= 2/3,~z_2= 0,~z_3= 1/3$, so it is not ideal.

 \end{example}
 \vspace{10pt}
 
 While the intersection of trees is not ideal or sharp, it still removes a significant number of fractional solutions from the linear relaxation compared to using formulations from \cite{veliborsPaper} or \cite{biggs2018optimizing}, leading to faster solve times as explored empirically in Section \ref{sec:numerical_experiments}. 

 \section{Strengthening formulations with binary split variables}
 
 We next present formulations that build upon the formulation from \cite{veliborsPaper}. In particular, these formulations use the binary variables from \cite{veliborsPaper}, which denote whether the feature vector is below each threshold in the tree. An advantage of this approach is its favorable branching behavior – setting a variable $x_{ij}=1$ will force all variables with a split threshold above this to also be 1, due to the ordering constraints $x_{ij}\leq x_{ij+1}$ (\ref{bigger_than_cosntraint}). In some cases, this results in a faster time to solve than the formulation in the previous section. We propose two ways to tighten this formulation to remove some of the fractional solutions, resulting in tighter linear relaxations and a faster time to solve in certain situations.

  \subsection{Tighter formulation from variable structure}
   \label{sec:tighter_variable}

   
 To tighten the formulation from \cite{veliborsPaper}, we exploit the  greater than or equal to representation of $\bm{x}$, which leads to larger groups of leaf variables being turned off when a split is made.  In \cite{veliborsPaper}, the $\bm{x}$ variables have consecutive 0's followed by consecutive 1's. In \cite{veliborsPaper}, if $x_{ij}=0$, this implies that all variables $z_l$ to the left of the split are equal to 0 (constraint \ref{right_constraint}). However, a stronger statement can be made. Due to the structure of $\bm{x}$, all variables with lower thresholds are also equal to 0, i.e., $x_{ik}=0 ~ \forall k < j$. This implies that variables $z_l$ to the left of splits with lower thresholds also must be equal to 0. 
 
\begin{figure}[]
  \centering
  \begin{subfigure}{0.48\textwidth}
    \centering
	\includegraphics[width=\textwidth]{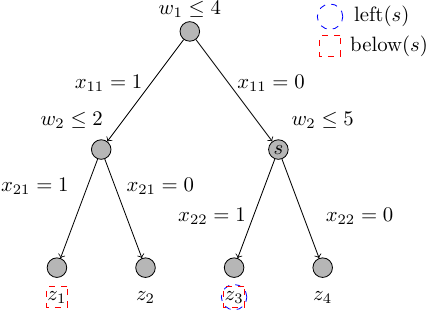}
	\caption{ The \textit{expset} formulation is tighter}
	\label{above_below_tree_example}
  \end{subfigure}
    \begin{subfigure}{0.48\textwidth}
    \centering
	\includegraphics[width=0.69\textwidth]{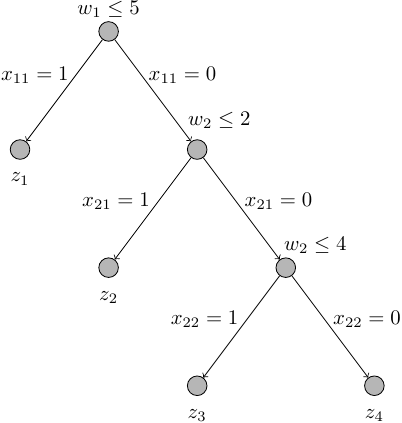}
	\caption{The \textit{elbow} formulation is tighter}
	\label{cuts_improve_example}
  \end{subfigure}
  \caption{Trees for which different formulations are tighter, and an illustration of the notation used in the \textit{expset} formulation, showing the sets $\textbf{left}(s)=\{3\}$ and $\textbf{below}(s)=\{1,3\}$. }
  \label{cut_above_below_figs}
\end{figure}

 As an illustrative example, we examine the tree in Figure \ref{above_below_tree_example}. If $w_2>5$ ($x_{22}=0$), then not only is the variable to the left of this split equal to 0, $z_3=0$, but also $z_1=0$ due to the constraint $x_{21} \leq x_{22}$ (constraint (\ref{bigger_than_cosntraint}) from \cite{veliborsPaper}). Rather than enforcing the relatively weak constraint from \cite{veliborsPaper}  that $z_3 \leq x_{22}$, it is tighter to directly enforce  $z_1 + z_3 \leq x_{22}$. Similarly, if $x_{ij}=1$, this implies that the variables $z_l$ to the right of any splits greater than the $j^{th}$ split are also set to 0. For example in Figure \ref{above_below_tree_example}, if $w_2 \leq 2$ ($x_{12}=1$), then not only is the variable to the right of this split equal to 0 ($z_2 = 0$), but also $z_4=0$, since the structure of $\bm{x}$ implies that $w_2 \leq 5$ ($x_{22}=1$).
 

 To formalize this logic, we introduce new sets $\textbf{below}(s)$ and $\textbf{above}(s)$. The set $\textbf{below}(s)$ contains all leaves to the left of splits with thresholds less than or equal to the threshold at split $s$ for a given tree. The set $\textbf{above}(s)$ contains all leaves to the right of leaves with a threshold greater than or equal to the threshold at split $s$. As such, for adjacent splits on the same feature, $s_{ij}$ and $s_{ij+1}$, we can define $ \textbf{below}(s_{ij+1})= \textbf{below}(s_{ij}) \cup \textbf{left}(s_{ij+1})$ and $ \textbf{above}(s_{ij})= \textbf{above}(s_{ij+1}) \cup \textbf{right}(s_{ij})$. For the smallest and largest splits, we have initial conditions that $\textbf{below}(s_{i1})=  \textbf{left}(s_{i1})$, and $\textbf{above}(s_{iK_i})=  \textbf{right}(s_{iK_i})$.  An equivalent pair of definitions are $\textbf{below}(s_{ij})= \bigcup_{k\leq j} \textbf{left}(s_{ik})$ and $\textbf{below}(s_{ij})= \bigcup_{k\geq j} \textbf{right}(s_{ik})$. An example of these sets is illustrated in Figure \ref{above_below_tree_example}. As a result, we can introduce a new formulation $Q^{expset}$, named after the notion of \textit{expanded sets}, by replacing (\ref{left_constraint}) and (\ref{right_constraint}) with the following constraints:
   \begin{subequations}
  \label{velibors_tighter_formulation}
\begin{align}
     Q^{expset}  = \{ \bm{x},y,\bm{z} ~ |~ &  \sum_{l \in \textbf{below}(s)} z_{l} \leq x_{V(s)C(s)}  \qquad \forall s ~\in  ~ \textbf{splits}(t)  \label{left_tight_constraint}\\
    & \sum_{l \in \textbf{above}(s)} z_{l} \leq 1- x_{V(s)C(s)}  \qquad \forall s ~\in  ~ \textbf{splits}(t) \label{right_tight_constraint} \\
    &  x_{ij} \leq x_{ij+1}  \qquad \forall i ~\in  ~ [p], ~ \forall j ~\in  ~ [K_i] \\
    & \sum_{l}^p z_{l} = 1, ~~~ y = \sum_{l = 1}^p s_l z_l  \\
    & \bm{x} \in [0,1]^{K_i} \qquad \forall i \in [d], ~\bm{z} \geq 0 \}.
\end{align}
  \end{subequations}
  
 Constraints (\ref{left_tight_constraint}) and (\ref{right_tight_constraint}) are the counterparts of (\ref{left_constraint}) and (\ref{right_constraint}). Constraint (\ref{left_tight_constraint}) enforces that when the condition at the split is not satisfied $x_{V(s)C(s)}=0$, the solution does not fall within a leaf to the left of any split in the tree with a lower threshold for the same feature, while constraint (\ref{right_tight_constraint}) enforces that all leaves to the right of greater splits are set to 0 if $x_{V(s)C(s)}=1$, as discussed previously. It can be shown that when intersected with a binary lattice on $\bm{x} \in \{0,1\}^p$, the feasible set of the MIO formulations (\ref{velibors_formulation}) and (\ref{velibors_tighter_formulation}) is the same. However, the linear relaxation, $Q^{expset}$ is generally a subset of $Q^{misic}$. This is shown in Proposition \ref{above_below_prop}, which formalizes the rationale given above.
\vspace{10pt}
\begin{proposition}
\label{above_below_prop}
The feasible sets associated with MIO formulations of $Q^{expset}$ and $Q^{misic}$ are equivalent, but the linear relaxation $Q^{expset}$ is a subset of $Q^{misic}$. Formally, $$Q^{expset} \cap (\{0,1\}^p  \times \mathbb{R}^{1+p})  = Q^{misic} \cap (\{0,1\}^p \times \mathbb{R}^{1+p}), ~ \text{but}~ Q^{expset} \subseteq Q^{misic}.$$
\end{proposition}


We provide a formal proof in Appendix \ref{proof_above_below_prop}. It can be shown that this formulation removes some fractional solutions from the LO relaxation of (\ref{velibors_formulation}). In particular, this will occur when there are multiple splits on the same feature within the tree. To illustrate this, suppose we have two splits on the same variable, $s$ and $s'$, where without loss of generality split $s'$ has the larger threshold. Define a reduced polyhedron that only includes the constraints related to these splits as follows:
\begin{align*}
&&\tilde{Q}^{expset}(s,s')= \{ \bm{x},\bm{z} ~ |&~ \sum_{l \in \textbf{below}(s)} z_{l} \leq x_{V(s)C(s)}, ~  \sum_{l \in \textbf{above}(s)} z_{l} \leq 1- x_{V(s)C(s)}, \\
&&&~\sum_{l \in \textbf{below}(s')} z_{l} \leq x_{V(s')C(s')}, ~  \sum_{l \in \textbf{above}(s')} z_{l} \leq 1- x_{V(s')C(s')},\\
&&& x_{V(s)C(s)} \leq  x_{V(s')C(s')}\},\\
&&\tilde{Q}^{misic}(s,s')= \{ \bm{x},\bm{z} ~ |&~ \sum_{l \in \textbf{left}(s)} z_{l} \leq x_{V(s)C(s)}, ~  \sum_{l \in \textbf{right}(s)} z_{l} \leq 1- x_{V(s)C(s)}, \\
&&&~\sum_{l \in \textbf{left}(s')} z_{l} \leq x_{V(s')C(s')}, ~  \sum_{l \in \textbf{right}(s')} z_{l} \leq 1- x_{V(s')C(s')}, \\
&&& x_{V(s)C(s)} \leq  x_{V(s')C(s')}\}.
\end{align*}

If we examine these polyhedrons, we see that the $\tilde{Q}^{expset}(s,s')$ is a strict subset of $\tilde{Q}^{misic}(s,s')$ when there are multiple splits on the same variable.
\vspace{10pt}

\begin{proposition}
\label{proof_tighter_subset}
Suppose we have two splits on the same variable, $s$ and $s'$, where $s'$ corresponds to the split with the larger threshold. Then  $$\tilde{Q}^{expset}(s,s') \subset \tilde{Q}^{misic}(s,s').$$
\end{proposition}

This is proved  in Appendix \ref{ap_proof_tighter_subset}. This proof involves exploring the potential relationships between splits $s$ and $s'$ (where split $s$ is a child of $s'$ in the tree, where $s'$ is a child of $s$, and where neither is a child of the other) and finding solutions $(\bm{x,z})$ that are in $\tilde{Q}^{misic}(s,s')$ but not in $\tilde{Q}^{expset}(s,s')$. An example that illustrates the strict subset is given in Example \ref{ex_tighter_stronger_than_cuts} from Section \ref{sec:comp_const}. In this example, we see that formulation (\ref{velibors_formulation}) has fractional solutions, while formulation (\ref{velibors_tighter_formulation}) has only integer solutions.

Generally, the more splits there are on the same feature in the tree, the more these constraints will tighten the formulation. At an extreme, we have the scenario where all splits in the tree are on the same feature. In the one-dimensional setting, it can be shown that the above formulation is ideal even for tree ensembles.

 
 

 \vspace{10pt}
 \begin{theorem} [Ideal formulation for one-dimensional tree ensembles]
 \label{ideal_one_dimension}
 The polyhedron defining a tree ensemble $\cap_{i=1}^T Q_i^{\text{expset}}$ is ideal if the feature is one-dimensional ($d=1$).
 \end{theorem}
 \vspace{10pt}
 
This result is proved in Appendix \ref{proof_ideal_one_dimension}. It follows by proving that the matrix representation of the polyhedron is totally unimodular. In particular, the matrix has a special structure whereby it is possible to provide a bi-coloring of the columns, such that the difference in row sums between the two groups is in $\{-1,0,1\}$. A result from \cite{ghouila1962caracterisation} proves that such a matrix is totally unimodular. A linear optimization formulation $\{\max \bm{c}'\bm{x} | A\bm{x}\leq \bm{b}\}$ has integer solutions if $b$ is integer and $A$ is a totally unimodular matrix \citep{schrijver1998theory}.
 

 The significance of Theorem \ref{ideal_one_dimension} is that it emphasizes the tightness of this formulation relative to other formulations that are not ideal in the one-dimensional scenario and have fractional solutions. In particular, in Example \ref{velibor_not_ideal}, we show that formulation from \cite{veliborsPaper} is not ideal in this case. In addition, the formulations from \cite{chen2021assortment} and \cite{kim2022reciprocity} do not have this property. 
 Furthermore, although this formulation isn't ideal when the input vector has multiple dimensions, we empirically show in Section \ref{sec:tighter_linear_relaxations} that the relaxation is tighter when the input vector is low dimensional.
  
 It is interesting to contrast this result with Theorem \ref{ideal_tree}. Theorem \ref{ideal_tree} states that the \textit{union of polyhedra} formulation is ideal for a single tree even with many features (there are similar results in \cite{chen2021assortment} and \cite{kim2022reciprocity} too). This contrasts with Theorem \ref{ideal_one_dimension}, which shows the \textit{expset} formulation is ideal for many trees but only if the ensemble has a single feature. This gives practitioners insight into the relative tightness of the different formulations. When there are many trees in the ensemble but relatively few variables, the \textit{expset} formulation is likely to be tighter. When there are few trees but many variables, the \textit{union of polyhedra} formulation is likely to be tighter. This formulation also provides an alternative way to strengthen the formulation from \cite{veliborsPaper} without introducing the large number of constraints that are introduced in \cite{kim2022reciprocity}, which lead to slow solve times in practice.

\blu
 Finally, we observe that this formulation can be extended when optimizing tree ensembles. In particular, $\textbf{below}(s)$ and $\textbf{above}(s)$ can be extended to include leaves from all trees that are below or above a particular split. Although this tightens the formulation further when optimizing tree ensembles, it presents practical implementation challenges. Since tree enseble data structures are typically arranged by tree (e.g., scikit-learn in Python), and adding constraints accross trees involves substantial restructuring of the data, or inefficient data access that can outweigh potential gains.
 \bla

 \subsection{Tighter formulation from nested branches}
\label{sec:elbow_cuts} 

 The relaxation of the formulation in the previous section still has some fractional extreme solutions, even when a single tree is being modeled over multiple features. These fractional extreme solutions often arise when there are nested splits, defined as follows: 
\vspace{10pt}
\begin{definition}
    A nested split occurs when a less-than (left) split is followed by a greater-than (right) split on the same feature on a path leading to a leaf or, alternatively, a less-than (left)  split follows a greater-than (right) split.
\end{definition}
\vspace{10pt}
 This is highlighted in the following example: 
 \vspace{10pt}
 \begin{example}[Nested splits that can be tightened]
 \label{elbow_cut_example}
 Consider a path to a leaf with nested splits shown in Figure \ref{frac_tree_example}. Suppose we model this using the formulation (\ref{velibors_formulation}) from \cite{veliborsPaper}: \[ \{x_1,x_2,z~| ~z \leq x_1,~ z \leq 1-x_2,~ x_2 \leq x_1,~  0 \leq  x_1, x_2 \leq 1, ~0 \leq z \}.\]

This has an extreme point $z=0.5, ~x_1=0.5,~ x_2=0.5$, as shown in Figure \ref{frac_ext_point_3d}. Consider the following reformulation: \[ \{x_1,x_2,z~| z \leq x_1-x_2, ~ 0 \leq  x_1,~ x_2 \leq 1, ~0 \leq z\}.\]
 
\begin{figure}[]
  \centering
  \begin{subfigure}{0.32\textwidth}
    \centering
	\includegraphics[width=0.75\textwidth]{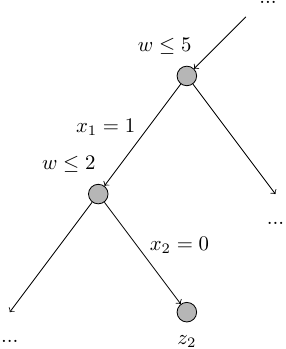}
	\caption{Tree segment}
	\label{frac_tree_example}
  \end{subfigure}
  \begin{subfigure}{0.32\textwidth}
    \centering
	\includegraphics[width=0.98\textwidth]{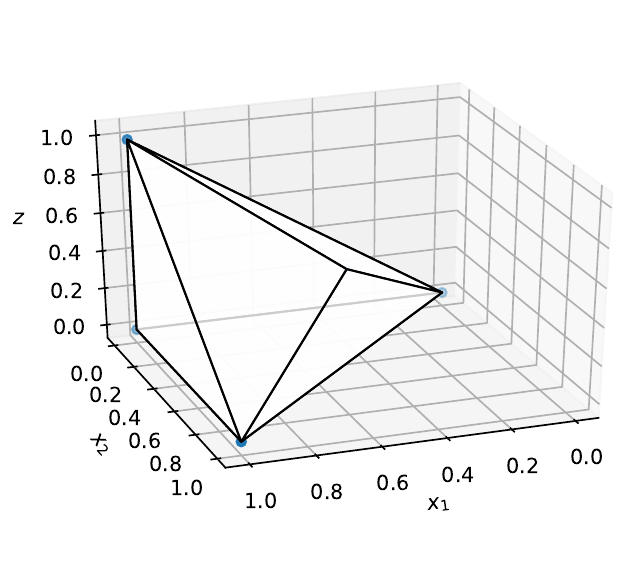}
\caption{$z \leq x_1,~ z \leq 1-x_2,~ x_2 \leq x_1$}
	\label{frac_ext_point_3d}
  \end{subfigure}
    \begin{subfigure}{0.32\textwidth}
    \centering
	\includegraphics[width=0.98\textwidth]{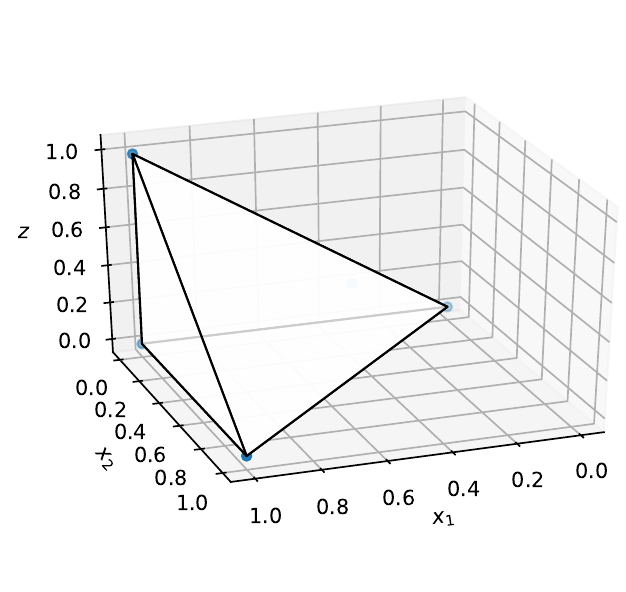}
\caption{$z \leq x_1-x_2$}
	\label{int_ext_point_3d}
  \end{subfigure}
  \caption{Example: cuts removing extreme point}
  \label{3d_plots}
\end{figure}

This is shown in Figure \ref{int_ext_point_3d}. As can be observed, this has removed the fractional extreme point, leaving only integer extreme points.  

 \end{example}
 \vspace{10pt}

More formally, we can characterize a valid set of constraints as follows: We define $\textbf{right\_parent}(s)$ as the set of splits that are above and to the right of split $s$ in the tree, with the additional requirement that these splits be on the same feature. That is, the split $s$ is a left child of another split on the same feature in the tree. For the splits in this set, the thresholds are necessarily larger.  We can also define $\textbf{left\_parent}(s)$ as the set of splits that are above and to the left of split $s$ for the same feature, for which the threshold is smaller. To illustrate this notation, in Figure \ref{cuts_improve_example} the split $w_2 \leq 2$ is the $\textbf{left\_parent}$ of the split $w_2 \leq 4$. We can generalize the constraints from Example \ref{elbow_cut_example} as follows:
\vspace{10pt}
\begin{definition} Nested split cuts:
       \begin{subequations}
  \label{tighter_formulation_cuts}
 \begin{align}
 & \sum_{l \in \textbf{right}(s)} z_{l} \leq x_{V(s')C(s')} - x_{V(s)C(s)}  ~~ \forall s ~\in  ~ \textbf{splits}(t),~ s' \in \textbf{right\_parent}(s), \label{left_tight_constraint_bend}\\
    &  \sum_{l \in \textbf{left}(s)} z_{l} \leq x_{V(s)C(s)} - x_{V(s')C(s')}  ~~ \forall s ~\in  ~ \textbf{splits}(t),~ s' \in \textbf{left\_parent}(s).\label{right_tight_constraint_bend} 
\end{align}
    \end{subequations}
\end{definition}

 If we define $Q^{elbow}$ as the polyhedron created by adding constraints (\ref{left_tight_constraint_bend}) and (\ref{right_tight_constraint_bend}) to formulation (\ref{velibors_formulation}) from \cite{veliborsPaper}, we can show that the relaxation of this formulation is tighter, while still having the same feasible region when $\bm{x}$ is restricted to a binary lattice, as shown in Proposition \ref{cuts_prop}.

 \vspace{10pt}
\begin{proposition}
\label{cuts_prop}
The feasible set associated with MIO formulations $Q^{elbow}$and $ Q^{misic}$ are equivalent, but linear relaxation $Q^{elbow}$ is a subset of $Q^{misic}$. Formally, $$Q^{elbow} \cap (\{0,1\}^p  \times \mathbb{R}^{1+p})  = Q^{misic} \cap (\{0,1\}^p \times \mathbb{R}^{1+p}), ~ \text{but}~ Q^{elbow} \subseteq Q^{misic}.$$
\end{proposition}

 This is proved formally in Appendix \ref{proof_cuts_prop}. As illustrated in Example \ref{elbow_cut_example}, the feasible region is often a strict subset when there are nested splits on the same feature ($Q^{elbow} \subset Q^{misic}$). This suggests that when there are more splits on the same features in the tree, there will be more of an improvement using the \textit{elbow} formulation over \cite{veliborsPaper}. This also often occurs if the tree has fewer features. This is explored empirically in Section \ref{sec:numerical_experiments}. However, simulation results suggest that the formulation is not ideal for tree ensembles with a single feature, unlike the \textit{expset} formulation.

 \subsection{Comparison of tightening constraints}
 \label{sec:comp_const}

In this section, we compare the relative tightness of the \textit{expset} and \textit{elbow} formulations ((\ref{velibors_tighter_formulation}) and (\ref{tighter_formulation_cuts}), respectively). We will show that when these constraints are added separately to formulation (\ref{velibors_formulation}) from \cite{veliborsPaper}, neither formulation is strictly tighter than the other. Rather, there are certain situations where one formulation is tighter than the other and vice versa, which we illustrate with examples.

A simple example where formulation (\ref{velibors_tighter_formulation}) is tighter than formulation (\ref{tighter_formulation_cuts}) is when there are multiple splits on the same variable, but they do not have a nested structure. For example, in the tree in Figure \ref{above_below_tree_example}, there are two splits on $w_2$, but these occur in different branches of the tree. In this situation, formulations (\ref{velibors_formulation}) and (\ref{tighter_formulation_cuts}) are the same since the constraints are added only for nested pairs of the same feature. Furthermore, formulation (\ref{tighter_formulation_cuts}) is not tight, but the formulation (\ref{velibors_tighter_formulation}) is tight.
\vspace{10pt}
\begin{example}[The \textit{expset} formulation is tighter than the \textit{elbow} formulation]
\label{ex_tighter_stronger_than_cuts}
For the tree given in Figure \ref{above_below_tree_example}, formulation (\ref{tighter_formulation_cuts}) (and formulation (\ref{velibors_formulation})) is:
\begin{align*}
\{\bm{x},\bm{z}~|~&x_{11} \geq z_1+z_2,
&& x_{21} \geq z_2,
&& x_{22}  \geq z_3, \\
& 1- x_{11} \geq z_3+z_4,
&& 1- x_{21} \geq z_2,
&& 1- x_{22} \geq z_4, \\
&x_{21}\leq x_{22},
&& z_1+z_2+z_3+z_4=1,
&& 0 \leq \bm{z}, ~0 \leq \bm{x} \leq 1\}.
\end{align*}
On the other hand formulation (\ref{velibors_tighter_formulation}) is: 
\begin{align*}
\{\bm{x},\bm{z}~|~&x_{11} \geq z_1+z_2,
&& x_{21} \geq z_2,
&& x_{22} \geq  \boxed{z_1}+ z_3, \\
& 1- x_{11} \geq z_3+z_4,
&& 1- x_{21} \geq z_2 +\boxed{z_4},
&& 1- x_{22} \geq z_4,\\
&x_{21}\leq x_{22},
&& z_1+z_2+z_3+z_4=1,
&&  0 \leq \bm{z}, 0 \leq \bm{x} \leq 1 \}.
\end{align*}
For convenience, the difference in the formulations has been highlighted. Formulation (\ref{tighter_formulation_cuts}) has fractional solutions $x_{11}=0.5, x_{21}=0.5, x_{22}=0.5, z_1=0,z_2=0.5,z_3=0,z_4=0.5$, and $x_{11}=0.5, x_{21}=0.5, x_{22}=0.5, z_1=0.5,z_2=0,z_3=0.5,z_4=0$, while formulation (\ref{velibors_tighter_formulation}) has only integer solutions.  The previous fractional solution violates the added constraints in formulation (\ref{velibors_tighter_formulation}). 
\end{example}
\vspace{10pt}
To further understand the difference between the constraints from formulations (\ref{tighter_formulation_cuts}) and (\ref{velibors_tighter_formulation}), it is useful to examine situations in which they are the same. In particular, suppose we have two nested splits on the same feature, such that  $s' \in \textbf{right\_parent}(s)$, as in the tree in Figure \ref{frac_tree_example}. We will examine constraints (\ref{left_tight_constraint}) and (\ref{right_tight_constraint}) and see when they imply the alternative constraint (\ref{left_tight_constraint_bend}). Specifically, we require that that $\textbf{above}(s)$ and $\textbf{below}(s') $ cover the whole set of 
leaves, that is, $\textbf{below}(s') \cup \textbf{above}(s) = p$. This is formally stated in Lemma \ref{theorem_cuts_above_below}.

\vspace{10pt}
\begin{lem}
\label{theorem_cuts_above_below}
    Suppose $s' \in \textup{\textbf{right\_parent}}(s)$. If $\textup{\textbf{below}}(s') \cup \textup{\textbf{above}}(s) = p$, 
    \begin{align*}
    Q^{misic} \bigcap  \sum_{l \in \textup{\textbf{below}}(s')} z_{l} \leq & x_{V(s')C(s')} \bigcap  \sum_{l \in \textup{\textbf{above}}(s)}  z_{l} \leq 1- x_{V(s)C(s)}  \\
    &\implies Q^{misic} \bigcap   \sum_{l \in \textup{\textbf{left}}(s)} z_{l} \leq x_{V(s)C(s)} - x_{V(s')C(s')}.
    \end{align*}
    Similarly, suppose $s' \in \textup{\textbf{left\_parent}}(s)$. If $\textup{\textbf{above}}(s') \cup \textup{\textbf{below}}(s) = p$, 
    \begin{align*}
    Q^{misic} \bigcap \sum_{l \in \textup{\textbf{below}}(s)} z_{l} \leq & x_{V(s)C(s)} \bigcap  \sum_{l \in \textup{\textbf{above}}(s')} z_{l} \leq 1- x_{V(s')C(s')}   \\ 
    &\implies Q^{misic} \bigcap \sum_{l \in \textup{\textbf{right}}(s)} z_{l} \leq x_{V(s)C(s)} - x_{V(s')C(s')}.
    \end{align*}
\end{lem}

This is proved in Appendix \ref{sec:proof_theorem_cuts_above_below}. The condition $\textbf{below}(s') \cup \textbf{above}(s) = p$ is satisfied when all splits above $s$ are on the same feature, or as an extreme case when the tree contains only one feature (the same condition as Theorem \ref{ideal_one_dimension}).  When these conditions are not met, including constraint (\ref{left_tight_constraint_bend}) will tighten the formulation. An example where this condition is not met and formulation (\ref{tighter_formulation_cuts}) is tighter than formulation (\ref{velibors_tighter_formulation}) occurs in Figure \ref{cuts_improve_example}. 


\vspace{10pt}
\begin{example}[\textit{Elbow} formulation is tighter than \textit{expset} formulation]

For the tree from Figure \ref{cuts_improve_example}, formulation (\ref{velibors_tighter_formulation}) is:
\begin{align*}
\{\bm{x},\bm{z}|~& x_{11} \geq z_1,
&& x_{21} \geq z_2,
&& x_{22} \geq z_2+ \boxed{z_3}, \\ 
& 1 - x_{11} \geq z_2+z_3+z_4,
&& 1- x_{21} \geq z_3+z_4,
&& 1- x_{22} \geq z_4,\\
& z_1+z_2+z_3+z_4=1,
&& x_{21}\leq x_{22},
&& 0 \leq \bm{x} \leq 1,~ 0 \leq \bm{z} \}.
\end{align*}

Formulation (\ref{tighter_formulation_cuts}) is:
\begin{align*}
\{\bm{x},\bm{z}|~& x_{11} \geq z_1,
&& x_{21} \geq z_2,
&& x_{22} \geq z_2, \\ 
& 1 - x_{11} \geq z_2+z_3+z_4,
&& 1- x_{21} \geq z_3+z_4,
&& 1- x_{22} \geq z_4,
~ ~ \boxed{x_{22}-x_{21} \geq z_3},  \\
&~ z_1+z_2+z_3+z_4=1, 
&& x_{21}\leq x_{22},
&& 0 \leq \bm{x} \leq 1, 
~ 0 \leq \bm{z} \}.
\end{align*}

For convenience, the difference in the formulations has been highlighted again. Formulation (\ref{velibors_tighter_formulation}) has a fractional solution  $x_{11}=0.5, x_{21}=0.5, x_{22}=0.5, z_1=0.5,z_2=0,z_3=0.5,z_4=0$, while formulation (\ref{tighter_formulation_cuts}) has only integer solutions. 

\end{example}
\vspace{10pt}

 Since each formulation has the advantage of removing different fractional solutions, including both sets of constraints can tighten the formulation further. We empirically explore how much these additional constraints tighten the LO relaxation for various datasets in Section \ref{sec:tighter_linear_relaxations}. 
 

\section{Numerical Experiments}
\label{sec:numerical_experiments}

\begin{table*}[]
\centering
\caption{Methods tested}
\begin{tabular}{lll}
Method & Reference & Comment  \\
\hline
\texttt{bigM} & \cite{biggs2018optimizing} & Not ideal for a single tree \\
\texttt{multilinear} & \cite{kim2022reciprocity} & Results in large number of constraints \\
\texttt{Mi\v{s}i\'{c}} & \cite{veliborsPaper}  – also (\ref{velibors_formulation}) & Not ideal for a single tree \\
\texttt{projected} & Formulation (\ref{projected_formulation}) & Ideal for a single tree \\
\texttt{expset} & Formulation (\ref{velibors_tighter_formulation}) & Ideal  for one feature tree ensembles \\ 
\texttt{elbow} & Formulation (\ref{velibors_formulation})+(\ref{tighter_formulation_cuts}) & Tighter than (\ref{velibors_formulation}) \\
\texttt{expset+elbow} & Formulation (\ref{velibors_tighter_formulation})+(\ref{tighter_formulation_cuts}) & Tightest formulation \\
\end{tabular}
\label{table:references}
\end{table*}

In this section, we study the numerical performance of the formulations on both simulated and real-world data. We study two scenarios of practical interest. The first involves the time taken to solve to optimality for an objective estimated by a tree ensemble. We then focus on finding tight \blu dual \bla bounds to this problem, obtained by solving the linear relaxation. 

\subsection{Experiments with tree ensembles}
\label{sec:tree_ensemble_experiments}

In this section, we examine the time taken to solve to optimality for a problem where the objective function is estimated using a random forest on simulated data. The random forest is trained on previous decisions where the reward is generated from a simple triangle-shaped function, where observed samples have added noise:

$$  r_i=\sum_{j=1}^d (1-|w_{ij}|) + d \cdot \epsilon_i.$$

For this problem, $r_i$ is a sampled reward, $w_{i} \sim U(-1,1)^d$ is a random decision vector with $d$ features, and $\epsilon_i \sim U(0,1)$ is added noise. There are no additional constraints placed on the variables other than those used to model the tree.  
We train a random forest from this data using \texttt{scikit-learn} \citep{scikit-learn}. The MIO formulations were solved using \texttt{Gurobi} solver, version 11.0.1 \citep{gurobi}, in \texttt{Python}, with a time limit of 30 minutes (1800s) for each trial but otherwise default parameters. The experiments were run on a MacBook Pro with an Intel 8-Core i9@2.4GHz with 32GB RAM.

\subsubsection{Small-scale forests}

We first show initial simulations on small-scale forests, to exhibit how the the full formulation from \cite{kim2022reciprocity}, denoted \texttt{multilinear}, is substantially slower than the other methods we tested. We then progress to larger forests, for which the \texttt{multilinear} cannot solve within the (1800s) limit for any instance. We compare formulation (\ref{projected_formulation}) denoted \texttt{projected} and formulation (\ref{tighter_formulation_cuts}) denoted \texttt{elbow}, to formulation (\ref{velibors_formulation}) from \cite{veliborsPaper}, denoted \texttt{Mi\v{s}i\'{c}}, and a formulation that uses the big-M method from \cite{biggs2018optimizing}, denoted \texttt{bigM}. This is summarized in Table \ref{table:references}.

\begin{table*}[]
\centering
\caption{Problem sizes for instance with 5 features, depth 8}
\begin{tabular}{lllll}
\hline
\# Trees            & Method & Constraints & Binary variables & Nonzeros \\
\hline
\multirow{4}{*}{1}  
  & multilinear     & 3712     & 189   & 119027     \\
  & projected    & 11       & 185   & 1623       \\
  & Mi\v{s}i\'{c}   & 1104     & 189   & 2341       \\
  & bigm    & 1104     & 553   & 2758       \\
  & elbow    & 585      & 189   & 2619       \\
  \hline
\multirow{4}{*}{2}  
  & multilinear     & 29277    & 379   & 1191149    \\
  & projected    & 22       & 376   & 3355       \\
  & Mi\v{s}i\'{c}   & 1119     & 374   & 4774       \\
  & bigm    & 2244     & 1124  & 5606       \\
  & elbow    & 1200     & 374   & 5266       \\
  \hline
\multirow{4}{*}{4}  
  & multilinear     & 223781   & 743   & 14246400   \\
  & projected    & 44       & 742   & 6606       \\
  & Mi\v{s}i\'{c}   & 2216     & 743  & 9433       \\
  & bigm    & 4428     & 11062 & 2218       \\
  & elbow    & 2381     & 743  & 10481      \\
  \hline
\multirow{4}{*}{8}  
  & multilinear     & 1839109  & 1511  & 211337514  \\
  & projected    & 88       & 1514  & 13588      \\
  & Mi\v{s}i\'{c}   & 4546     & 1506  & 19326      \\
  & bigm    & 9036     & 4526  & 22574      \\
  & elbow    & 4888     & 1506  & 21358      \\
\end{tabular}
\label{table:problem_size_small}
\end{table*}

For the small-scale simulations, we calculate the solve time to optimality with an increasing number of trees in the forest and an increasing number of features. We increase the number of trees according to  $\{1,2,4,8\}$. We use default parameters and a maximum depth for each tree of 8. For these parameters, each tree has an average of $189$ leaves. We show problem sizes of the formulations when there are 5 features in Table \ref{table:problem_size_small}. This shows the number of constraints,  binary variables, and the sparsity of the constraint matrix with the number of nonzero entries.  We note the very large number of constraints for \texttt{multilinear}, even for relatively small problem instances. With this many constraints, even creating the model in \texttt{Gurobi} can be prohibitively slow. As mentioned earlier, the number of constraints in \texttt{projected} formulation is substantially smaller, while the number of binary variables is also less than the other formulations.

\begin{figure}
  \centering
	\includegraphics[width=0.52\textwidth]{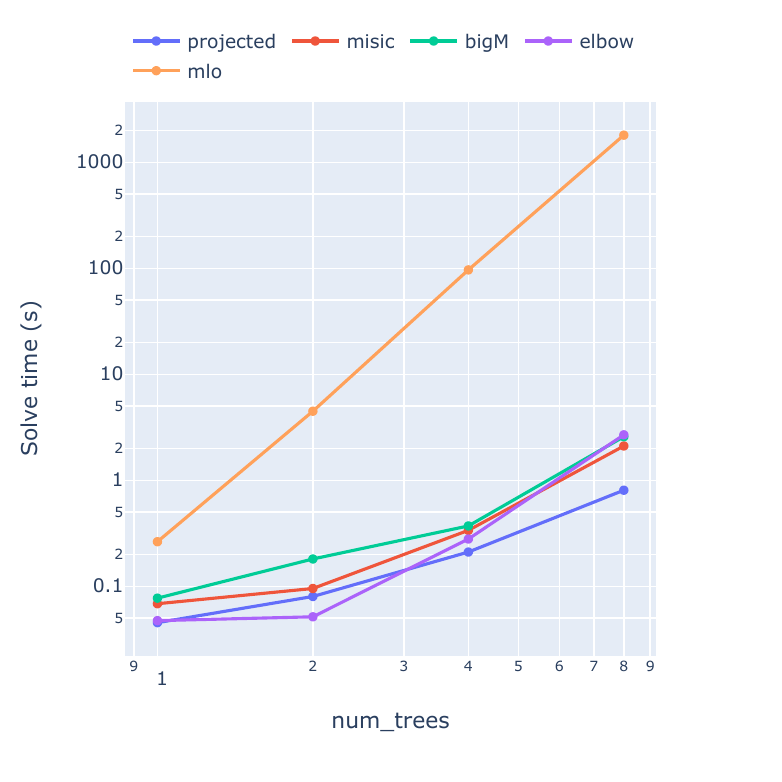}
	\caption{Solve time for small-scale forest}
\label{fig:small_scale_time_taken}
  \end{figure}

The time taken to solve to optimality is given in Figure \ref{fig:small_scale_time_taken}, on a log-log axis for clarity with \texttt{multilinear} referenced as \texttt{mlo}. We highlight that problems that can be solved in approximately 1 second using the other methods cannot be solved within the 30-minute time window for the \texttt{multilinear} formulation. However, we recognize that the solve time might be improved by employing additional techniques such as delayed constraint generation or removal of redundant constraints, which we don't implement for consistency with the other formulations and since this adds substantially to the complexity of the implementation.  

\subsubsection{Larger forests}
For the experiments on larger forests, we increased the depth from 8 to 20, increasing the number of leaves per tree from 189 to 2761 on average. We additionally study the effect of varying the number of features from 1 to 5. We also increase the number of trees according to  $\{1, 2, 4, 8, 16, 32\}$. We show the size of problem instances in terms of variables and constraints in Appendix \ref{sec:problem_size}. We repeat the experiment for 10 randomly generated datasets for each forest size and number of features.

\begin{figure}
  \centering
  \begin{subfigure}{0.42\textwidth}
    \centering
	\includegraphics[width=\textwidth]{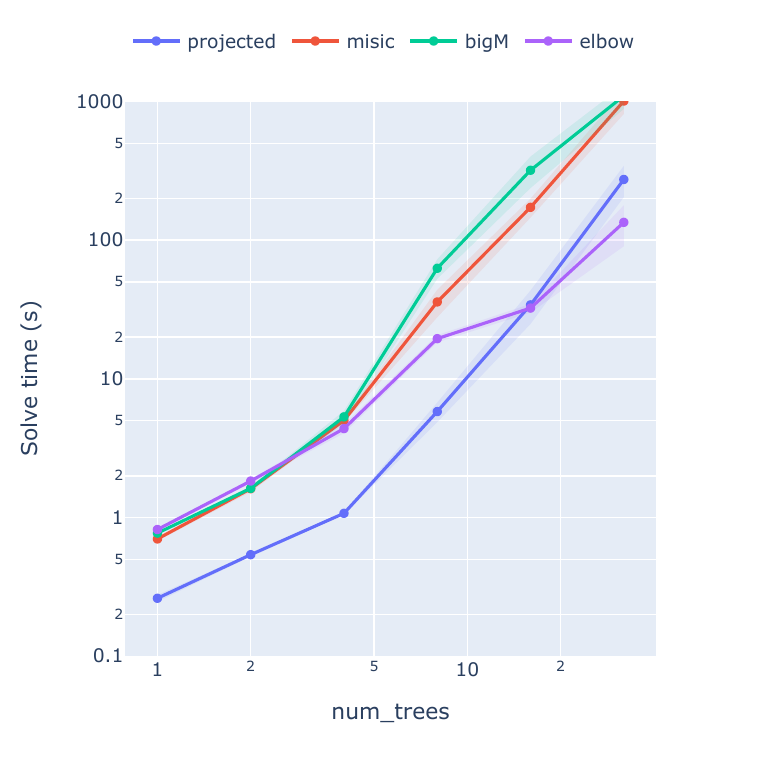}
	\caption{1 feature}
  \end{subfigure}
    \begin{subfigure}{0.42\textwidth}
    \centering
	\includegraphics[width=\textwidth]{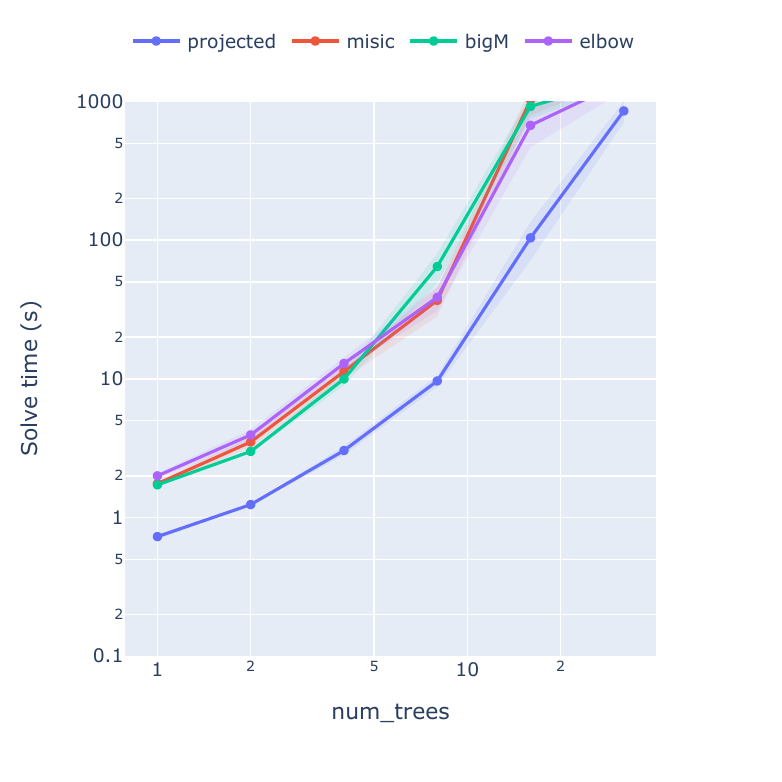}
	\caption{2 features}
  \end{subfigure}
    \begin{subfigure}{0.42\textwidth}
    \centering
	\includegraphics[width=\textwidth]{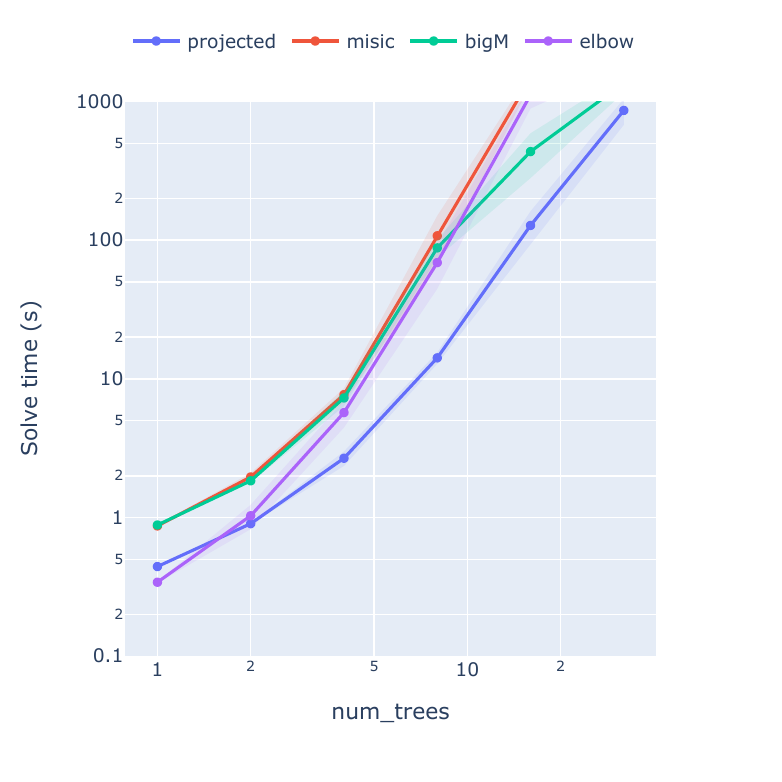}
	\caption{3 features}
  \end{subfigure}
     \begin{subfigure}{0.42\textwidth}
    \centering
	\includegraphics[width=\textwidth]{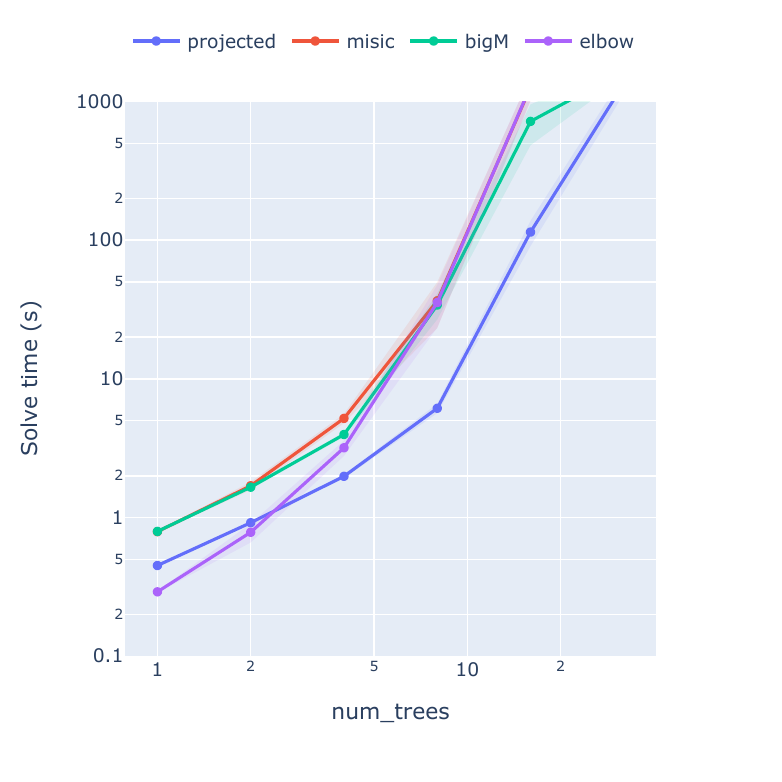}
	\caption{4 features}
  \end{subfigure}
     \begin{subfigure}{0.42\textwidth}
    \centering
	\includegraphics[width=\textwidth]{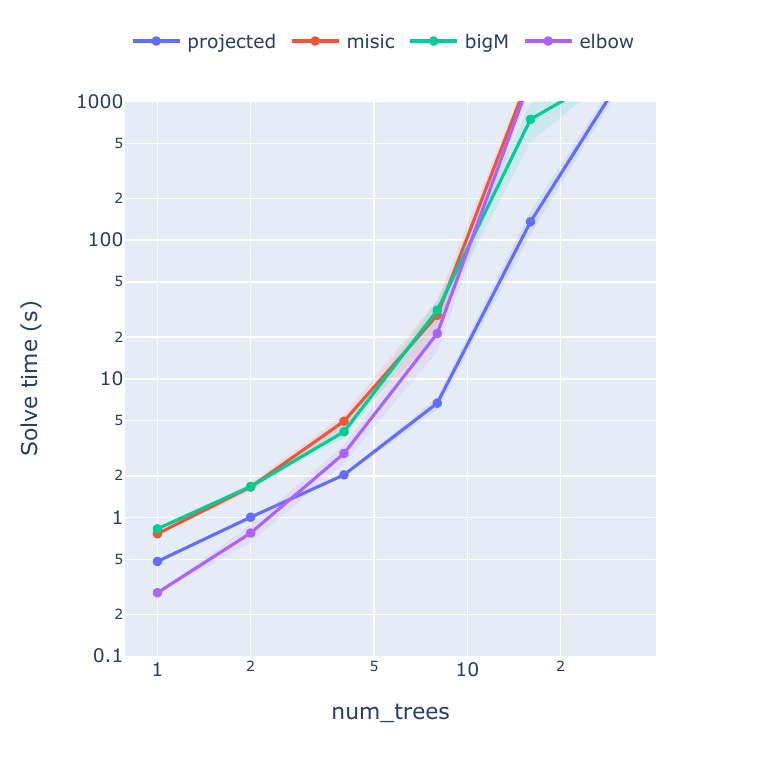}
	\caption{5 features}
  \end{subfigure}
  \caption{Time taken to solve to optimality for random forests of varying sizes, depth 20}
  \label{fig:time_taken_rf}
\end{figure}

In Table \ref{table:time_taken_rf} we observe the time taken to solve optimally for different-sized trees. Each result is averaged over 50 trials: 10 trials for each input vector of 1 to 5 dimensions. We note that the average time taken includes instances that didn't reach optimality, recorded as the maximum time allocated (1800s), so it is in fact a truncated mean. The percentage of instances that didn't reach optimality is recorded in the last four columns. As can be seen, the projected formulation is on average three to four times faster, and it finds an optimal solution more often within the given time. 


\begin{table*}
    \centering
    \caption{Time taken to solve to optimality, depth 20}
    \begin{tabular}{cccccccc}
       \hline
       & & \multicolumn{6}{c}{ \# Trees }      \\
        \cmidrule(lr){3-8} 
       &   & 1 & 2 & 4 & 8 & 16 & 32  \\
       \hline
       \multirow{4}{*}{Time taken (s)}  & \texttt{projected} & 0.47  &  0.92 & 2.16 &  8.50 & 103.30  & 983.29 \\
         & \texttt{Mi\v{s}i\'{c}} & 0.98 &  2.09 & 6.83 & 49.14 & 1111.25 & 1552.09 \\
         &  \texttt{bigM}  & 1.00  & 1.96 & 6.15 & 56.16 & 628.49 & 1477.53\\
         & \texttt{elbow} & 0.75 & 1.67 & 5.82 & 36.82 & 914.28  & 1363.65 \\
       \hline
         \multirow{4}{*}{\% greater 1800s} & \texttt{projected} & 0 & 0 & 0 & 0 & 0 & 32  \\
         & \texttt{Mi\v{s}i\'{c}} & 0 & 0 & 0 & 0 & 42 & 76\\
         & \texttt{bigM} & 0 & 0 & 0 & 0 & 14 & 66\\
         & \texttt{elbow} & 0 & 0 & 0 & 0 & 38 & 70\\
    \end{tabular}
\label{table:time_taken_rf}
\end{table*}

Figure  \ref{fig:time_taken_rf} shows the results further broken down by the number of features, plotted on a log-log axis for clarity. We observe that the \texttt{elbow} formulation is often faster for tree ensembles with few trees. This might be useful in applications where many MIO problems need to be solved rapidly, such as policy iteration in reinforcement learning with tree-based value function approximations. We also observe a substantial solve time improvement using the \texttt{elbow} formulation when there is one feature, which agrees with the results presented in Section \ref{sec:elbow_cuts}. 


\subsubsection{Tighter linear relaxations}
\label{sec:tighter_linear_relaxations}
A problem of practical interest is finding tight \blu dual \bla bounds for optimization problems with an objective estimated by a tree ensemble. For large problem instances, finding an optimal solution can be prohibitively slow, considering that MIO formulations often exhibit exponential solve times. The relative quality of a fast heuristic solution can be assessed if an \blu upper/lower \bla bound on the objective can be found when \blu maximizing/minimizing. \bla Another application of dual bounds is the verification of the robustness of a machine learning model \citep{carlini2017towards, dvijotham2018training}, whereby an optimization problem is solved over local inputs to find maximally different output. Since finding the exact worst case change can be prohibitively slow for large instances, a bound is often used instead.

\begin{figure}
  \centering
  \begin{subfigure}{0.42\textwidth}
    \centering
	\includegraphics[width=\textwidth]{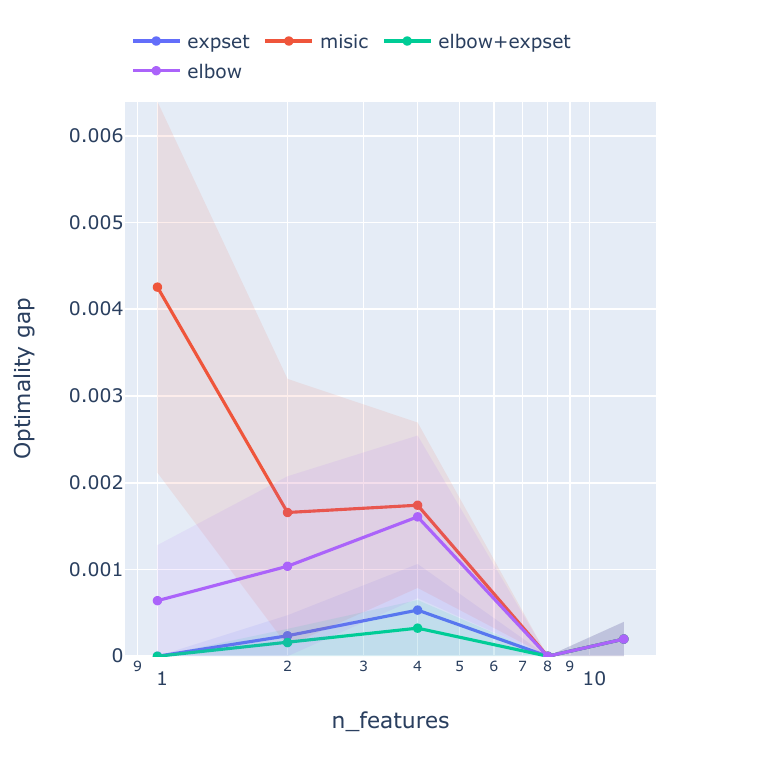}
	\caption{2 trees}
  \end{subfigure}
    \begin{subfigure}{0.42\textwidth}
    \centering
	\includegraphics[width=\textwidth]{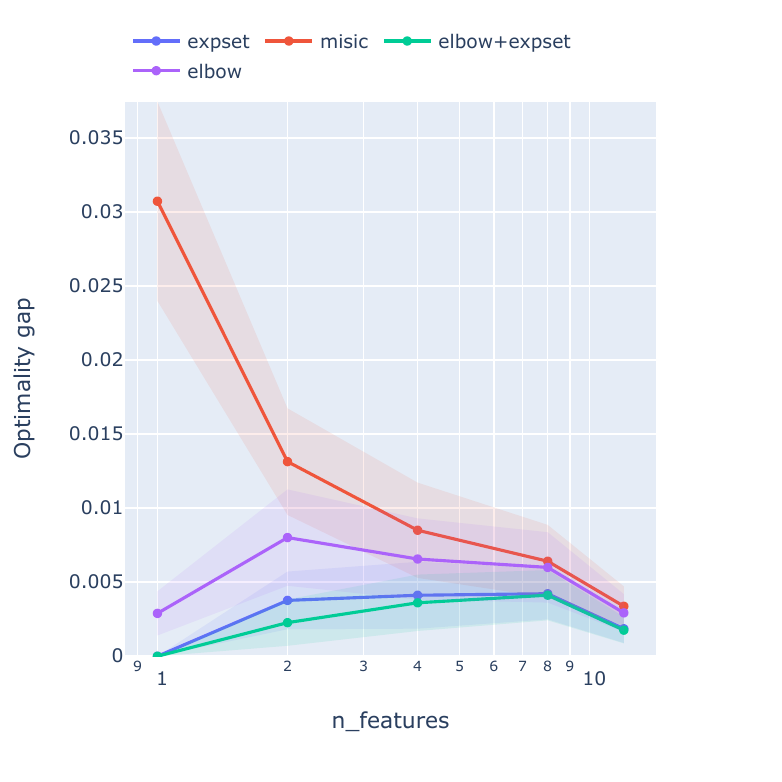}
	\caption{4 trees}
  \end{subfigure}
     \begin{subfigure}{0.42\textwidth}
    \centering
	\includegraphics[width=\textwidth]{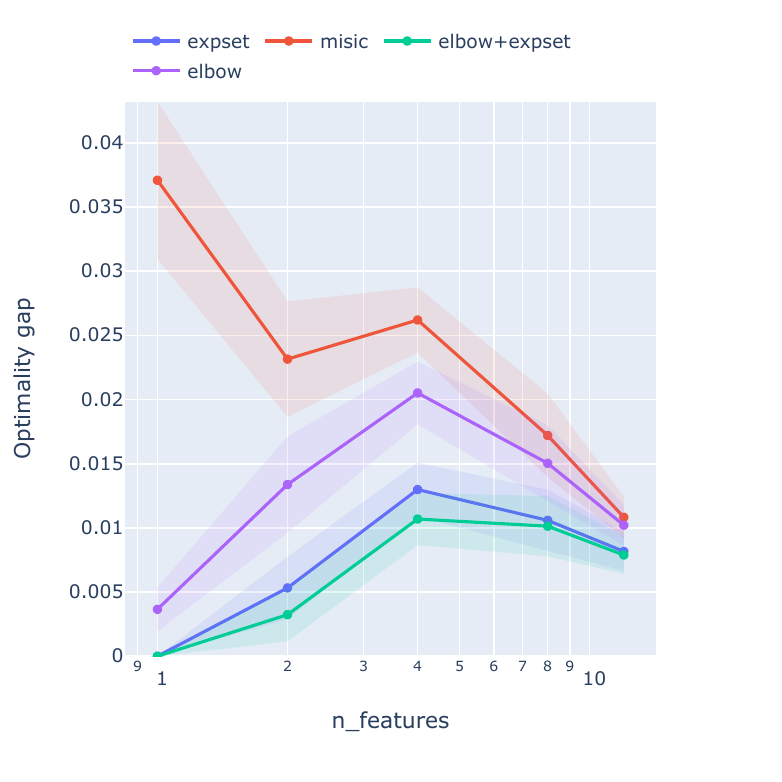}
	\caption{6 trees}
  \end{subfigure}
    \begin{subfigure}{0.42\textwidth}
    \centering
	\includegraphics[width=\textwidth]{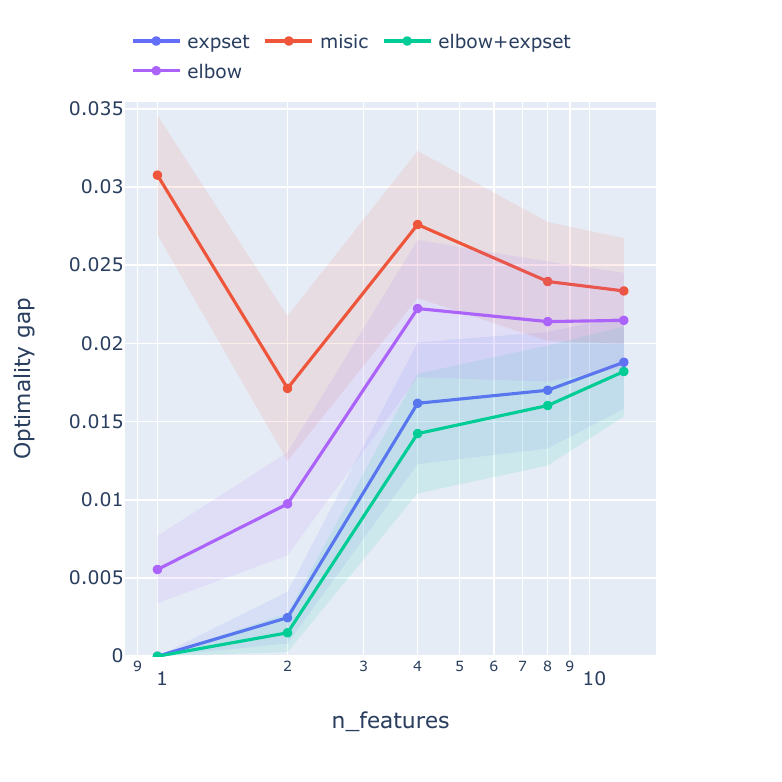}
	\caption{8 trees}
  \end{subfigure}
      \begin{subfigure}{0.42\textwidth}
    \centering
	\includegraphics[width=\textwidth]{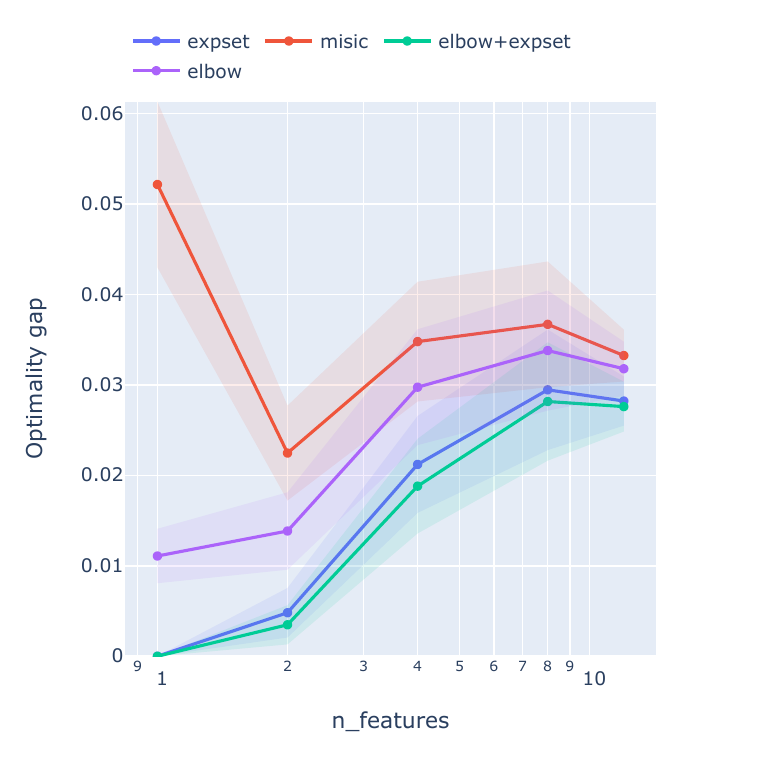}
	\caption{10 trees}
  \end{subfigure}
  \caption{Tightness of linear relaxation}
  \label{fig:tighter_linear_relaxation}
\end{figure}

We analyze the formulations from Section \ref{sec:tighter_variable} by analyzing the tightness of the linear relaxation. We compare formulations that use the same variables, specifically formulation (\ref{velibors_tighter_formulation},  \texttt{expset}), formulation (\ref{velibors_formulation}, \texttt{Mi\v{s}i\'{c}}), and   (\ref{tighter_formulation_cuts}, \texttt{elbow}). Additionally, we test a formulation that has both of the tightening constraints (\texttt{expset+elbow}). 
We use the same data-generating process as in Section \ref{sec:tree_ensemble_experiments}, except rather than solving to find an optimal integer solution, we solve only the linear relaxation. For these experiments, we use forests with $\{2,4,6,8,10\}$ trees, and increase the features according to $\{1,2,4,8,12\}$. Again, we repeat each experiment with 10 randomly generated datasets.

Figure \ref{fig:tighter_linear_relaxation} shows the optimality gap \blu fraction\bla, calculated from the difference between the objective of the linear relaxation and an optimal value, as the number of features increases. We observe the effect of Theorem \ref{ideal_one_dimension}, whereby for tree ensembles with one feature, formulations based on \texttt{expset} are ideal. Moreover, for problems with relatively few features, the formulation is significantly tighter than formulation \texttt{Mi\v{s}i\'{c}}, whereas when the number of features is larger, the improvement is smaller. This is likely due to more features being associated with fewer splits per feature. 
We note that in isolation, the constraints introduced in \texttt{expset} have a greater effect in tightening the formulation than those introduced in \texttt{elbow}, although combining both results in the tightest formulations. We also empirically observe that the \texttt{elbow} formulation is not ideal even in the single feature case.

\blu An alternative to examining the linear relaxation would be to restrict Gurobi to solving the root node only. This would allow Gurobi to use its presolve techniques and generate cuts to improve the solution. We have focused on the linear relaxation to isolate the tightness of the respective formulations without the effect being complicated by additional factors. Furthermore, even the root node can take a considerable amount of time to solve for some problem instances. \bla


\subsection{Real-world data}


\begin{figure}
  \centering
  \begin{subfigure}{0.42\textwidth}
    \centering
	\includegraphics[width=\textwidth]{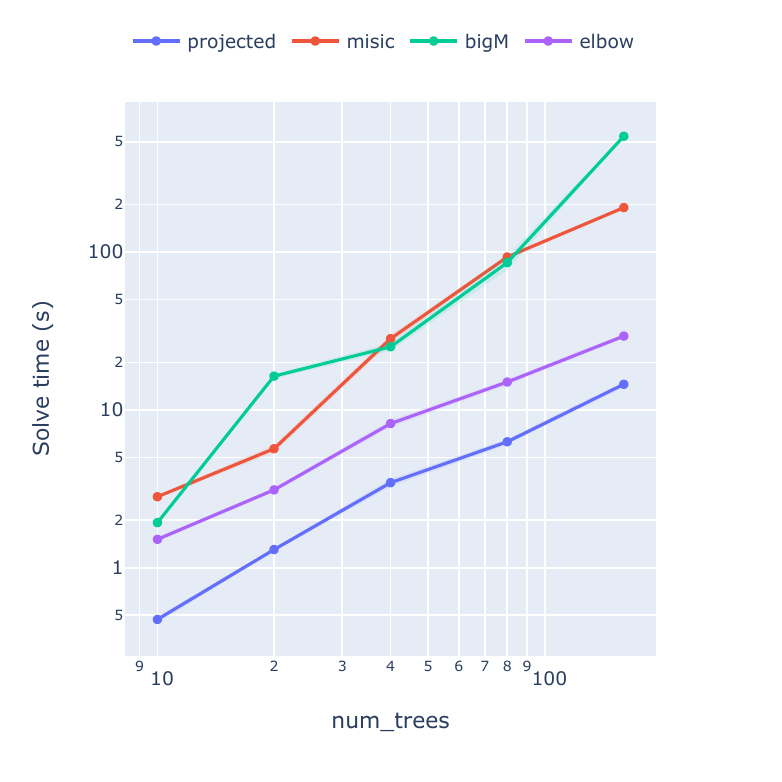}
	\caption{1 feature}
 \label{fig:concrete_1_feat}
  \end{subfigure}
    \begin{subfigure}{0.42\textwidth}
    \centering
	\includegraphics[width=\textwidth]{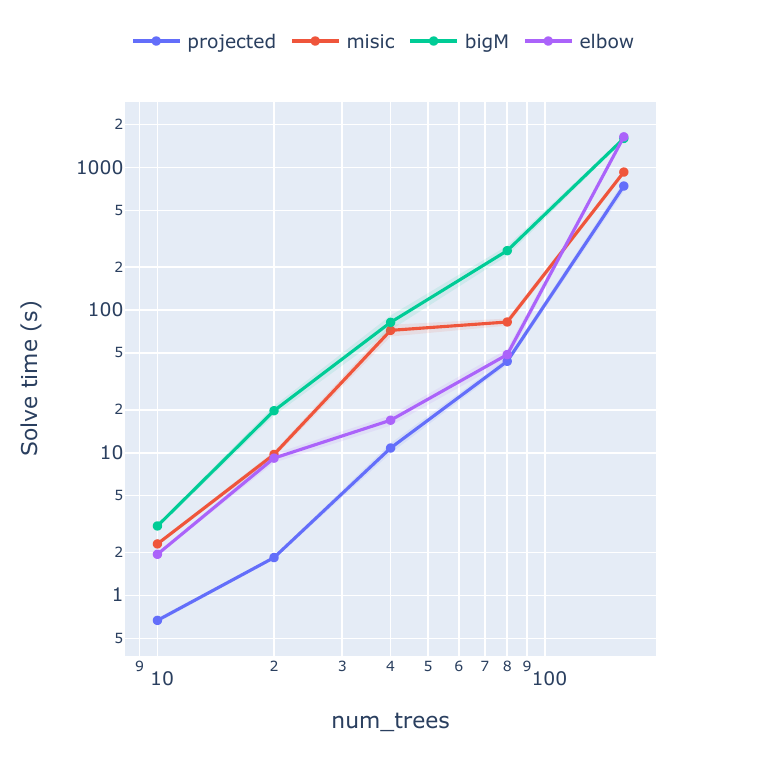}
	\caption{3 features}
  \end{subfigure}
    \begin{subfigure}{0.42\textwidth}
    \centering
	\includegraphics[width=\textwidth]{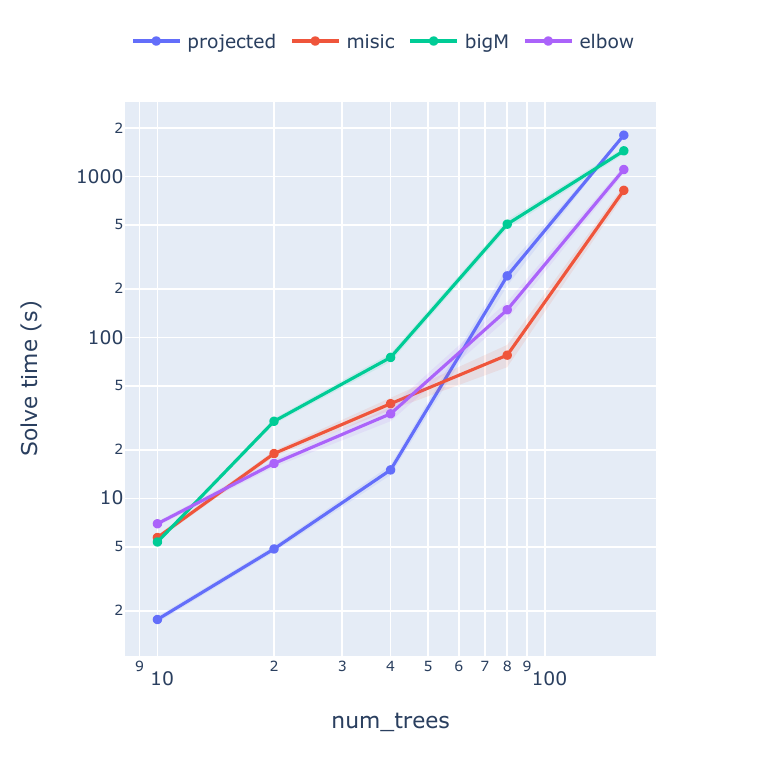}
	\caption{5 features}
  \label{fig:concrete_5_feat}
  \end{subfigure}
     \begin{subfigure}{0.42\textwidth}
    \centering
	\includegraphics[width=\textwidth]{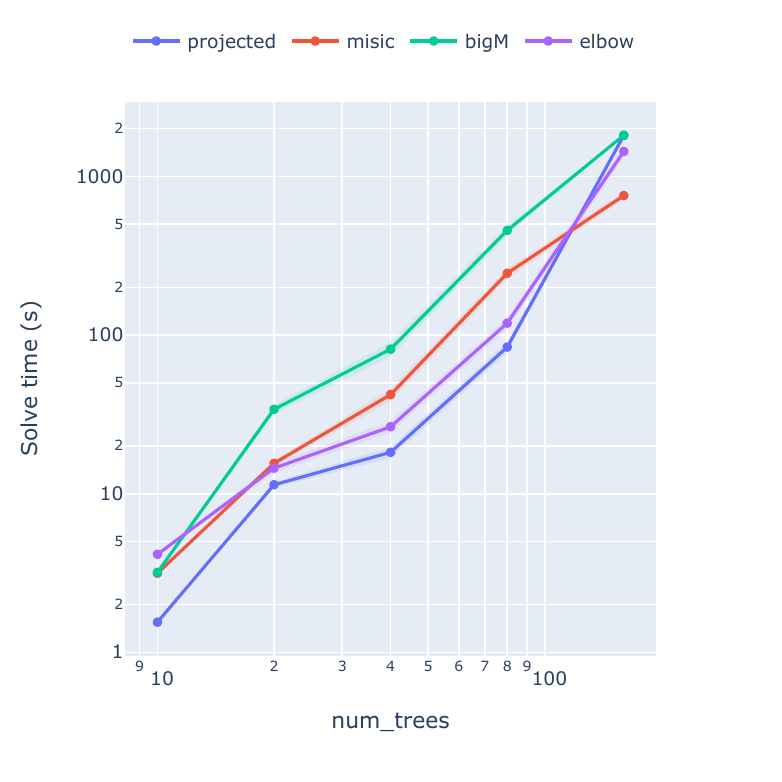}
	\caption{7 features}
  \end{subfigure}
  \caption{Time taken to solve to optimality for random forests of varying sizes \texttt{concrete} data}
  \label{fig:time_taken_rf_concrete}
\end{figure}

\begin{figure}
  \centering
  \begin{subfigure}{0.42\textwidth}
    \centering
	\includegraphics[width=\textwidth]{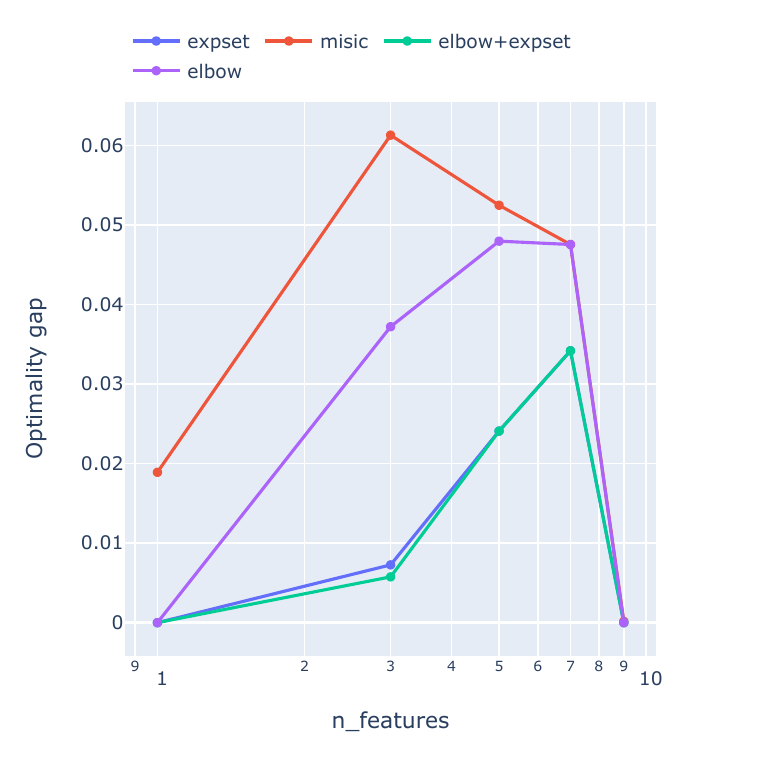}
	\caption{10 trees}
  \end{subfigure}
    \begin{subfigure}{0.42\textwidth}
    \centering
	\includegraphics[width=\textwidth]{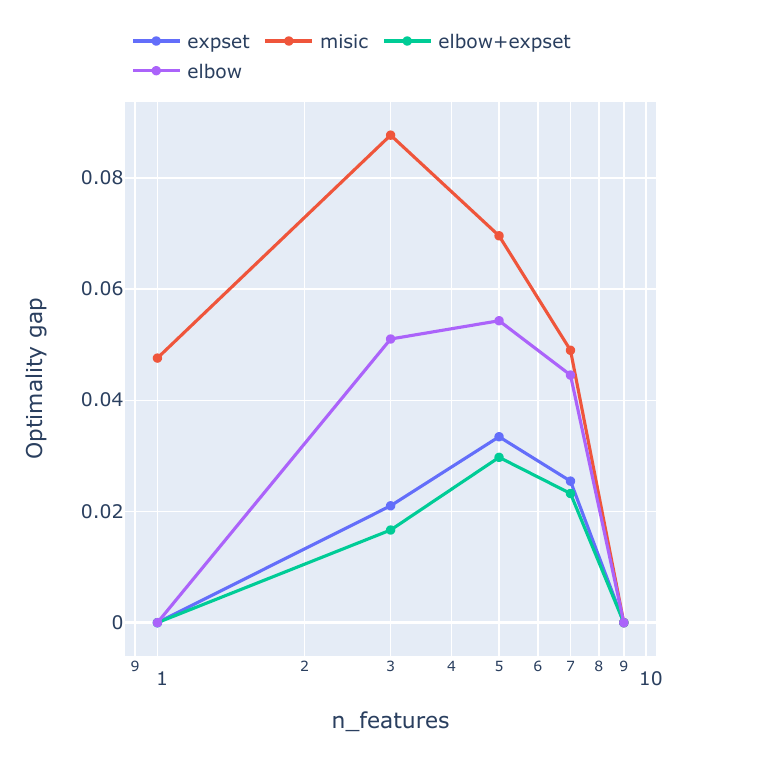}
	\caption{20 trees}
  \end{subfigure}
    \begin{subfigure}{0.42\textwidth}
    \centering
	\includegraphics[width=\textwidth]{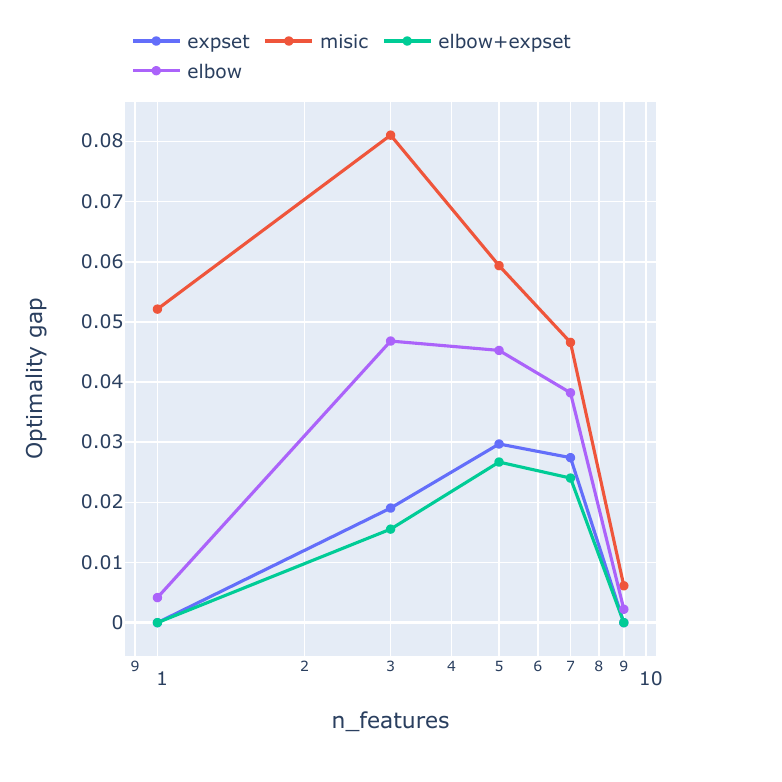}
	\caption{40 trees}
  \end{subfigure}
     \begin{subfigure}{0.42\textwidth}
    \centering
	\includegraphics[width=\textwidth]{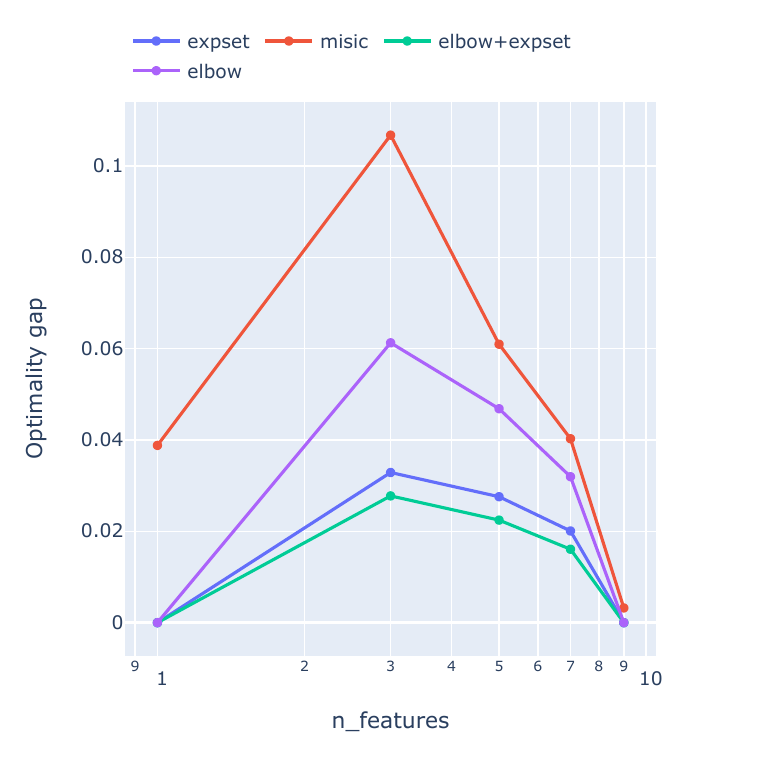}
	\caption{80 trees}
  \end{subfigure}
  \caption{Tightness of linear relaxation for random forests of varying sizes \texttt{concrete} data}
  \label{fig:relaxation_concrete}
\end{figure}


\begin{figure}
  \centering
  \begin{subfigure}{0.42\textwidth}
    \centering
	\includegraphics[width=\textwidth]{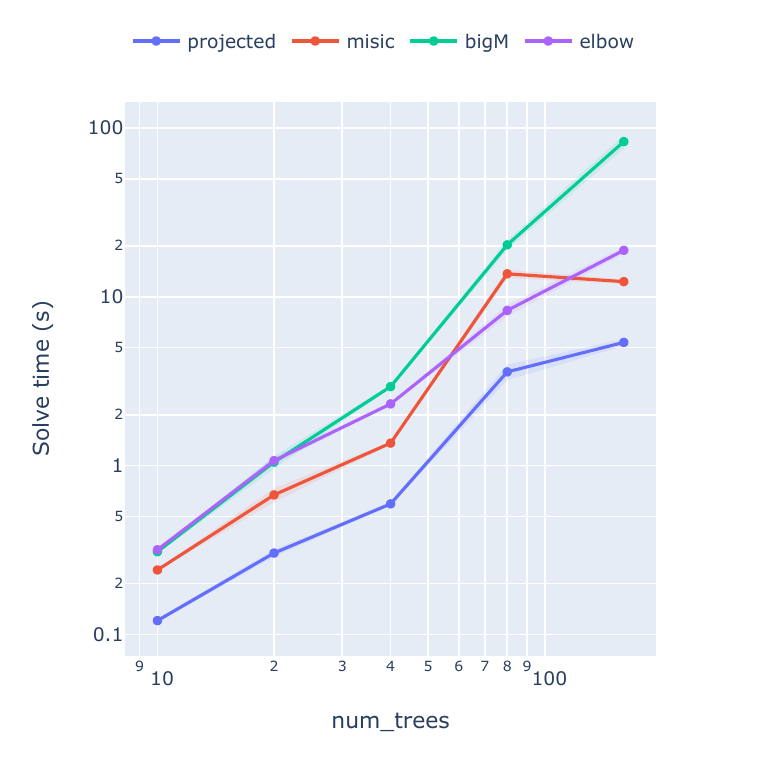}
	\caption{1 feature}
  \end{subfigure}
    \begin{subfigure}{0.42\textwidth}
    \centering
	\includegraphics[width=\textwidth]{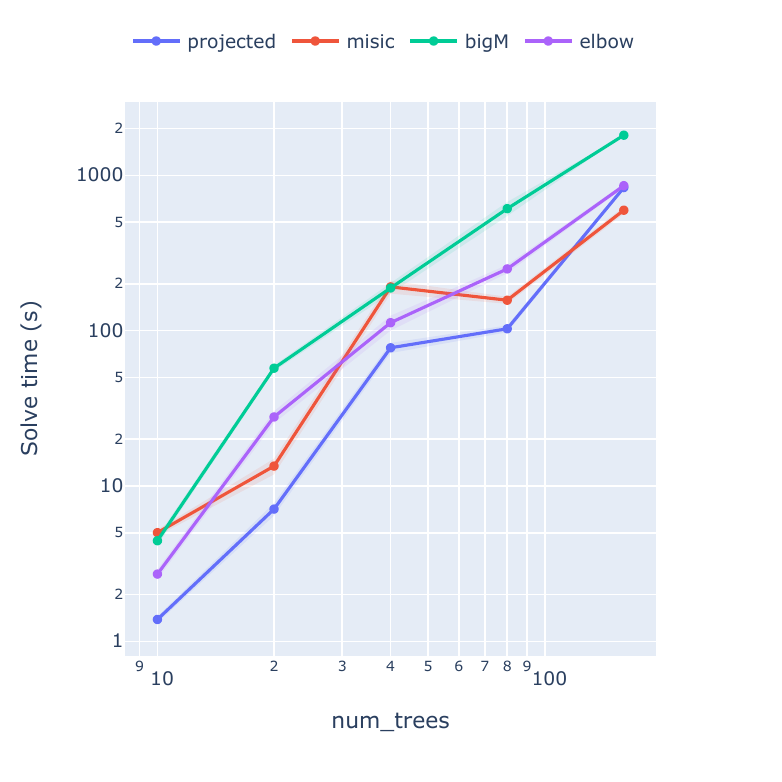}
	\caption{5 features}
  \end{subfigure}
     \begin{subfigure}{0.42\textwidth}
    \centering
	\includegraphics[width=\textwidth]{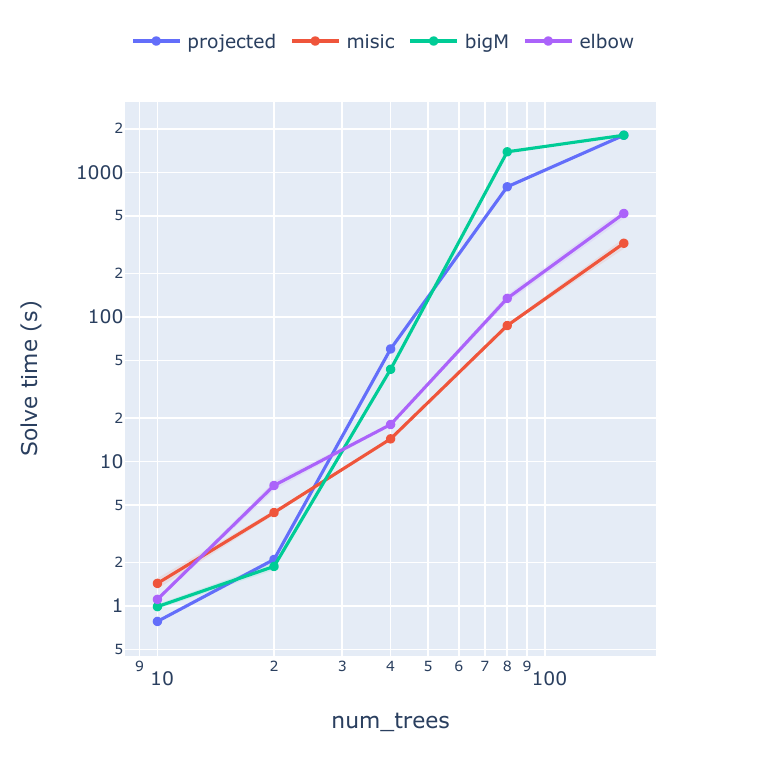}
	\caption{10 features}
 \label{fig:wine_10_feat}
  \end{subfigure}
  \caption{Time taken to solve to optimality for random forests of varying sizes \texttt{winequalityred} data}
  \label{fig:time_taken_rf_wine}
\end{figure}

\begin{figure}
  \centering
    \begin{subfigure}{0.42\textwidth}
    \centering
	\includegraphics[width=\textwidth]{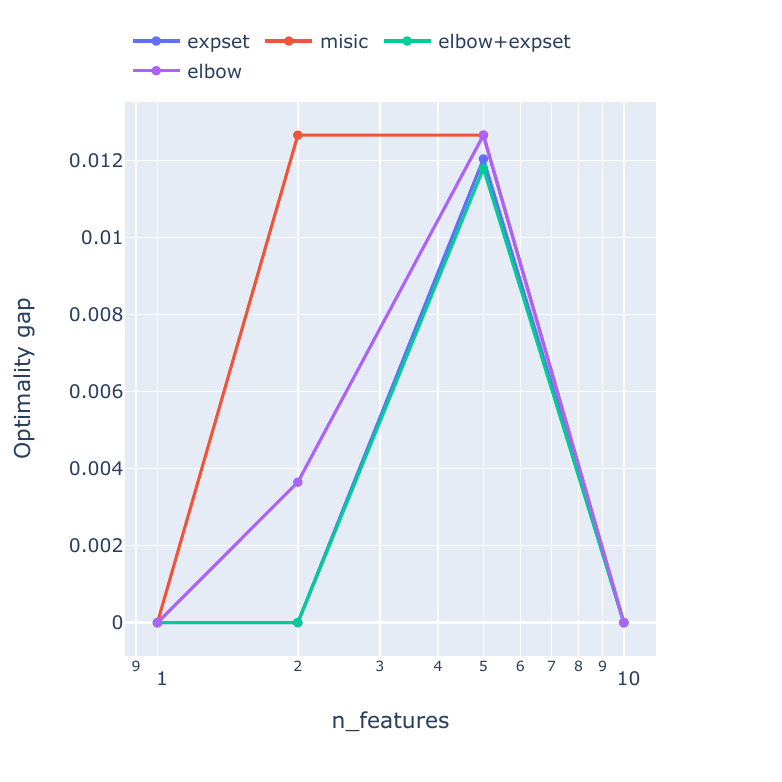}
	\caption{20 trees}
  \end{subfigure}
     \begin{subfigure}{0.42\textwidth}
    \centering
	\includegraphics[width=\textwidth]{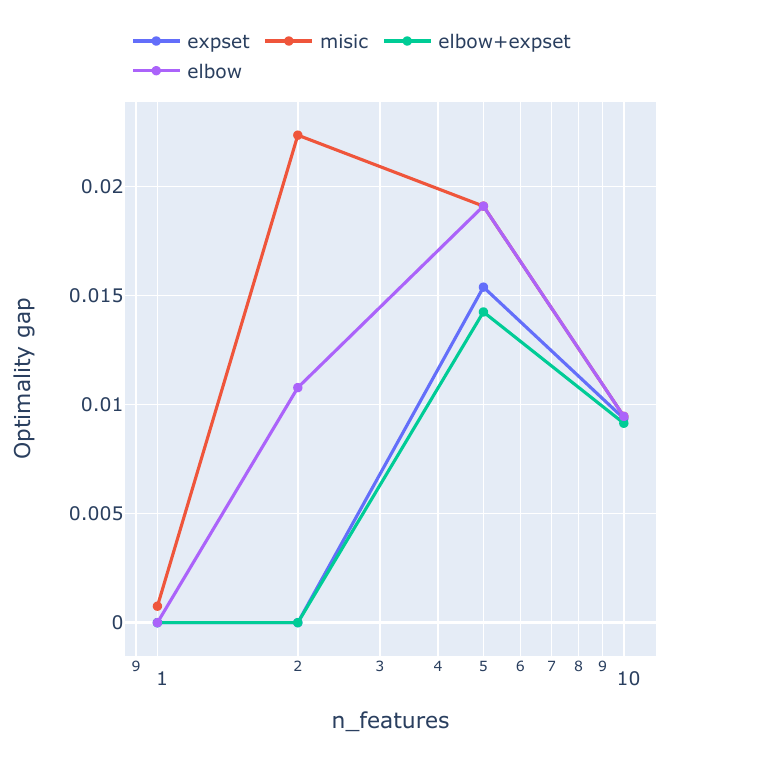}
	\caption{40 trees}
  \end{subfigure}
     \begin{subfigure}{0.42\textwidth}
    \centering
	\includegraphics[width=\textwidth]{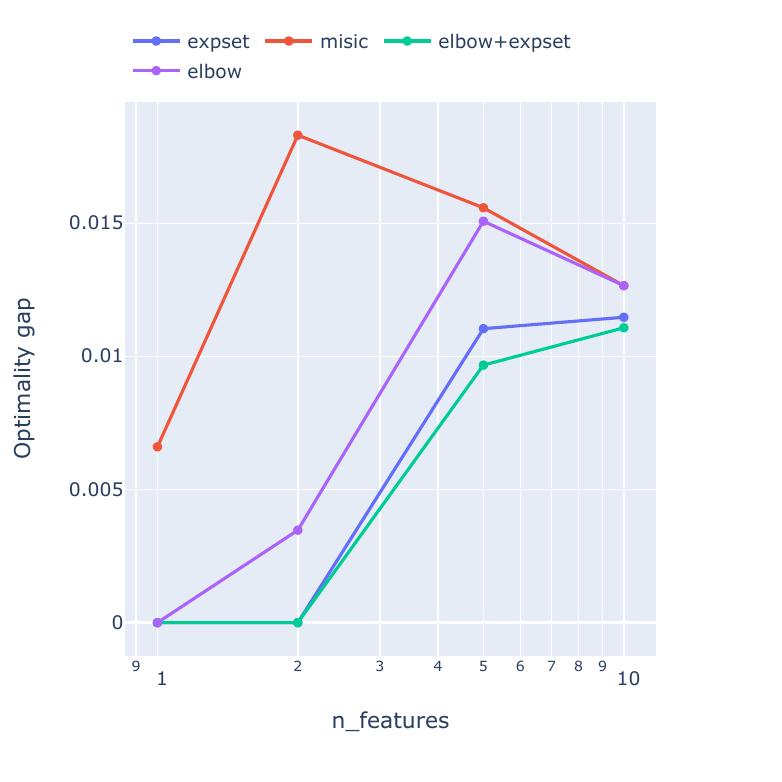}
	\caption{80 trees}
  \end{subfigure}
  \caption{Tightness of linear relaxation for random forests of varying sizes \texttt{winequalityred} data}
  \label{fig:relaxation_rf_wine}
\end{figure}

We also study some datasets used to benchmark tree ensemble solve times used in \cite{veliborsPaper}. In particular, we study the \texttt{concrete} dataset \citep{yeh1998modeling}, with 1030 observations. The dependent variable is the compressive strength of concrete, with independent variables being the characteristics of the concrete mix. \footnote{Cement, BlastFurnaceSlag, FlyAsh, Water, Superplasticizer, CoarseAggregate, FineAggregate, Age.} Optimization aims to find the concrete with the highest compressive strength. We also study the \texttt{winequalityred} dataset \cite{cortez2009modeling}, with 1599 observations. The dependent variable is the quality of the wine, while the independent variables are characteristics of the wine. \footnote{fixed acidity, volatile acidity,	citric acid,	residual sugar,	chlorides, free sulfur dioxide, total sulfur dioxide,	density, pH, sulfates, alcohol.} As such, the optimization problem is to choose characteristics of the wine such that the quality is maximized.  

\subsubsection{Solve time}
\label{sec:solve_time_real_data}
We explore the solve time for different formulations of different size random forest tree ensembles $\{10, 20, 40, 80, 160\}$ and varying feature vector dimension $\{1,3,5,7\}$ for \texttt{concrete} and $\{1,5,10\}$ for \texttt{winequalityred}. To test the effect of dimension, we use the first $k$ features to predict the output. As in the previous section, we set the maximum solve time to be 30 minutes (1800s).

The results for \texttt{concrete} and  \texttt{winequalityred} are in Figures \ref{fig:time_taken_rf_concrete} and \ref{fig:time_taken_rf_wine}, respectively. We observe that for both datasets, the \texttt{projected} formulation performs relatively better than the formulation from \cite{veliborsPaper} for  instances where the feature vector has a lower dimension (fewer features). On the other hand, for instances with a larger number of features, the formulation  \cite{veliborsPaper} can be faster to solve. 
Furthermore, the projected formulation (\ref{projected_formulation}) appears to be relatively faster for formulations with a small number of trees, which is particularly pronounced in Figures  \ref{fig:concrete_5_feat} and \ref{fig:wine_10_feat}. This is potentially an extension of Theorem \ref{ideal_tree}; if (\ref{projected_formulation}) is ideal for a single tree, it is also potentially relatively tighter for a small number of trees. Again, this might have applications where many smaller problems need to be solved quickly, such as in reinforcement learning. For these datasets, the performance of the \texttt{elbow} formulation is generally comparable to \cite{veliborsPaper}, although there are improvements in the \texttt{concrete} dataset when there are few features.

\subsubsection{Tightness of linear relaxation}

For the \texttt{concrete} and \texttt{winequalityred} datasets, we also compare the tightness of the linear relaxations for the \texttt{concrete} and \texttt{winequalityred} datasets in Figures \ref{fig:relaxation_rf_wine} and \ref{fig:relaxation_concrete}. Across both datasets, we observe a similar outcome to the synthetic data experiments, whereby \texttt{elbow+expset} is generally the tightest, followed by \texttt{expset}, and finally, the original \texttt{Mi\v{s}i\'{c}} formulation. We also observe that generally, the difference diminishes when there are more features in the data, potentially because there are fewer splits per feature, which is typically where the new formulations remove fractional points.

\section{Conclusions and future work}

In this paper, we have proposed a variety of new mixed-integer optimization formulations for modeling the relationship between an input feature vector and the predicted output of a trained decision tree. We have introduced formulations that build on the variable structure from \cite{veliborsPaper} and formulations that use the input feature directly. We have shown that these formulations are provably tighter than existing formulations in certain scenarios. We have shown conditions where these formulations are ideal, which gives further practical insight into when different formulations might be advantageous, depending on the number of trees in the ensemble and the number of features the problem has. In addition to these theoretical insights, we have given experimental conditions where the different formulations succeed both in terms of the time taken to solve to optimality and the tightness of the corresponding linear relaxations. While the experimental results do not always fully agree with the theoretical findings or intuition due to the complex operations of commercial MIO solvers, we have identified situations where each different formulation has advantages and laid the groundwork for future computational studies.

For future work, an interesting avenue is exploring the relationship between the formulations we provide and different polyhedral constraints. While, in general, the formulations we provide are not ideal when combined with additional constraints, there may be special cases when they are, or at least cuts that can be introduced to remove some of the fractional solutions.  
An additional promising direction is exploring formulations that encode leaf selection using a logarithmic number of binary variables, potentially improving computational performance by significantly reducing the number of integer variables, though this approach may introduce additional complexity in mapping constraints to the binary encoding.
Another relevant extension would be to explore MIO formulations that can handle non-rectangular decision trees, such as oblique or hyperplane-based splits. Our formulations exploit the axis-aligned nature of the leaf partitions and can not be applied directly. It is unclear whether formulations tighter than those presented in \cite{biggs2018optimizing} exist for this setting.

\backmatter








\section*{Declarations}

\subsection*{Competing interests}
The authors have no competing interests to declare that are relevant to the content of this article.

The code for the numerical studies can be found in \url{https://anonymous.4open.science/r/Tightness-of-Prescriptive-Tree-Based-Mixed-Integer-Optimization-4D13}









\bibliography{sn-bibliography}

\begin{appendices}
\input{appendix} 





\end{appendices}


\end{document}

%% file: appendix.tex
\section{Proof Theorem \ref{ideal_tree}}
\label{proof_ideal}
\proof{}

We prove this by applying Fourier-Motzkin elimination to formulation (\ref{union_polyhedra_formulation}) to eliminate all $\bar{\bm{w}}_{l}$, and showing that we arrive at formulation (\ref{projected_formulation}). An overview of the technique can be found in \cite{hooker2011logic}. For convenience, recall $Q^{ext}$:
 \begin{subequations}
\begin{align}
Q^{\text{ext}} = \{  \bm{w}, y , \bar{\bm{w}}, \bar{\bm{y}}, \bm{z} ~| ~ & u_{li} z_l \geq \bar{w}_{li} \qquad \forall i \in [d], ~\forall l \in [p] \label{union_upper_bound_app}\\
&  b_{li} z_l \leq \bar{w}_{li}  \qquad \forall i \in [d],  ~\forall l \in [p] \label{union_lower_bound_app}\\
& \bar{y}_l = s_l z_l,\qquad \forall l \in [p] \label{union_score_app} \\
& \sum_{l = 1}^p z_l =1,  \label{union_choose_one_set_app}\\
& w_i = \sum_{l = 1}^p \bar{w}_{li}  \qquad \forall i \in [d] \label{union_form_w_app} \\ 
& y = \sum_{l = 1}^p \bar{y}_l  \label{union_form_y_app} \\
& z_l \in [0,1]  \qquad \forall l \in [p] \label{union_z_app} \}.
\end{align}
 \end{subequations}

To eliminate $\bar{\bm{w}}_{l}$, we will use induction. To be more precise, we will show how to eliminate $\bar{w}_{1i},...,\bar{w}_{pi}$ for a single feature $i$, but applying the same procedure to the other features is identical. For notational brevity, let us define $Q^{const}$ as the set of constraints that do not feature $\bar{w}_{1i},...,\bar{w}_{pi}$ and do not change with elimination.
\begin{subequations}
\begin{align}
Q^{const} = \{  \bm{w}, y , \bar{\bm{w}}, \bar{\bm{y}}, \bm{z} ~~|~  &  u_{lj} z_l \geq \bar{w}_{lj} \qquad \forall j \neq i  , ~\forall l \in \{1,...,p\}\\
&  b_{lj} z_l \leq \bar{w}_{lj}  \qquad \forall j \neq i ,~ \forall l \in \{1,...,p\} \\
& \bar{y}_l = s_l z_l,\qquad \forall l \in \{1,...,p\} \\
& y = \sum_{l = 1}^p \bar{y}_l  \\
& \sum_{l = 1}^p z_l =1,  \\
& z_l \in \{0,1\}  \qquad \forall l \in \{1,...,p\}  \}.
\end{align}
\end{subequations}

Define $Q^{proj}_k$ as the polyhedron resulting from applying Fourier-Motzkin elimination $k$ times on $Q^{ext}$ to eliminate $\bar{w}_{1i},...,\bar{w}_{ki}$. We propose $Q^{proj}_k$ is
\begin{subequations}
\begin{align}
Q^{proj}_k = Q^{const} \cap   \{  \bm{w}, y , \bar{\bm{w}}, \bar{\bm{y}}, \bm{z} ~~|~ & \sum_{l =1}^{k} b_{ki} z_l \leq w_i - \sum_{l \in \{k+1,...,p\}} \bar{w}_{li}  \\
  & \sum_{l =1}^{k}  u_{ki} z_l \geq w_i - \sum_{l \in \{k+1,...,p\}} \bar{w}_{li}   \\
&  u_{li} z_l \geq \bar{w}_{li}  \qquad~\forall l \in \{k+1,...,p\}\\
&  b_{li} z_l \leq \bar{w}_{li}  \qquad ~ \forall l \in \{k+1,...,p\} \}.
\end{align}
\end{subequations}

As the inductive step, if we apply Fourier-Motzkin elimination to $Q^{proj}_k$ to eliminate $\bar{w}_{(k+1)i}$, we will show that $Q^{proj}_{k+1}$ is the resulting polyhedron. First, we establish the base case, that applying Fourier-Motzkin elimination on  $Q^{ext}$ to eliminate $\bar{w}_{1i}$ results in $Q^{proj}_1$. 

To apply Fourier-Motzkin elimination, we rearrange all constraints involving $\bar{w}_{1i}$ into greater than constraints $\bar{w}_{1i} \geq G_j(\bm{w}, y , \bar{\bm{w}}, \bar{\bm{y}}, \bm{z})$ or less than constraints  $\bar{w}_{1i} \leq L_{j'}(\bm{w}, y , \bar{\bm{w}}, \bar{\bm{y}}, \bm{z})$. We eliminate these constraints and replace them with $L_{j'}(\bm{w}, y , \bar{\bm{w}}, \bar{\bm{y}}, \bm{z}) \geq  G_j(\bm{w}, y , \bar{\bm{w}}, \bar{\bm{y}}, \bm{z})$ for all combinations $j$ and $j'$. As a result, the new constraints formed are
 \begin{subequations}
\begin{align}
 &&b_{1i} z_1 \leq \bar{w}_{1i} , ~ \bar{w}_{1i} \leq u_{1i} z_1   &&\implies  & b_{1i} z_1 \leq  u_{1i} z_1  \label{eq:redundant} \\
&&    w_i - \sum_{l \in \{2,...,p\}} \bar{w}_{li} \leq  \bar{w}_{l1}, ~  \bar{w}_{1i}  \leq u_{1i} z_1 && \implies &  u_{1i} z_1 \geq w_i - \sum_{l \in \{2,...,p\}} \bar{w}_{li}  \label{eq:upper_temp}  \\
&&  w_i  - \sum_{l \in \{2,...,p\}} \bar{w}_{li} \geq  \bar{w}_{l1} , ~ \bar{w}_{1i}  \geq b_{1i} z_1 && \implies &  b_{1i} z_1 \leq w_i - \sum_{l \in \{2,...,p\}} \bar{w}_{li} \label{eq:lower_temp} \\
&&   w_i - \sum_{l \in \{2,...,p\}} \bar{w}_{li} \leq  \bar{w}_{l1}, ~  \bar{w}_{l1} \leq  w_i  - \sum_{l \in \{2,...,p\}} \bar{w}_{li},  \nonumber  \\ 
&& \implies  w_i - \sum_{l \in \{2,...,p\}} \bar{w}_{li} && \leq & ~~w_i - \sum_{l \in \{2,...,p\}} \bar{w}_{li}, \label{eq:redundant_also}
\end{align}
 \end{subequations}
 
\noindent where the constraint (\ref{eq:redundant}) is formed by combining (\ref{union_upper_bound_app}) and (\ref{union_lower_bound_app}), constraint (\ref{eq:upper_temp}) is from (\ref{union_upper_bound_app}) and (\ref{union_form_w_app}),  (\ref{eq:lower_temp}) is from (\ref{union_lower_bound_app}) and (\ref{union_form_w_app}), and (\ref{eq:redundant_also}) is from (\ref{union_form_w_app}). By definition, constraint (\ref{eq:redundant}) is redundant and can be eliminated,  since $b_{1i} \leq u_{1i}$, as can (\ref{eq:redundant_also}). As a result, the polyhedra is:
\begin{subequations}
\begin{align}
Q^{proj}_1 = Q^{const} \cap \{  \bm{w}, y , \bar{\bm{w}}, \bar{\bm{y}}, \bm{z} ~~|~   & b_{1i} z_1 \leq w_i - \sum_{l \in \{2,...,p\}} \bar{w}_{li}  \\
  & u_{1i} z_1 \geq w_i - \sum_{l \in \{2,...,p\}} \bar{w}_{li}   \\
&  u_{li} z_l \geq \bar{w}_{li}  \qquad~\forall l \in \{2,...,p\}\\
&  b_{li} z_l \leq \bar{w}_{li}  \qquad ~ \forall l \in \{2,...,p\}   \}.
\end{align}
\end{subequations}

We can apply the same logic to prove the inductive step. If we apply Fourier-Motzkin elimination to $Q^{proj}_k$ to eliminate $\bar{w}_{(k+1)i}$ we get
 \begin{subequations}
\begin{align}
b_{(k+1)i} z_{k+1} &\leq \bar{w}_{(k+1)i},  \bar{w}_{(k+1)i}  \leq  u_{(k+1)i} z_{k+1}  \implies b_{(k+1)i} z_{k+1}  \leq  u_{(k+1)i} z_{k+1} \label{eq:redundant_2} \\
b_{(k+1)i} z_{k+1} & \leq \bar{w}_{(k+1)i}, ~\bar{w}_{(k+1)i}    \leq w_i - \sum_{l \in \{k+2,...,p\}} \bar{w}_{li} - \sum_{l =1}^{k} b_{ki} z_l \nonumber \\
& \implies   b_{(k+1)i} z_{k+1}   \leq w_i - \sum_{l \in \{k+2,...,p\}} \bar{w}_{li} - \sum_{l =1}^{k} b_{ki} z_l  \label{eq:lower_temp_2}\\
 u_{(k+1)i} z_{k+1} & \geq \bar{w}_{(k+1)i}, ~\bar{w}_{(k+1)i}  \geq w_i - \sum_{l \in \{k+1,...,p\}} \bar{w}_{li} -\sum_{l =1}^{k}  u_{ki} z_l  \nonumber \\
 & \implies     u_{(k+1)i} z_{k+1}    \geq w_i - \sum_{l \in \{k+1,...,p\}} \bar{w}_{li} -\sum_{l =1}^{k}  u_{ki} z_l \label{eq:upper_temp_2} \\
  w_i - &\sum_{l \in \{k+2,...,p\}} \bar{w}_{li} - \sum_{l =1}^{k} b_{ki} z_l \geq  \bar{w}_{(k+1)i}    , ~\bar{w}_{(k+1)i}  \geq w_i - \sum_{l \in \{k+1,...,p\}} \bar{w}_{li} -\sum_{l =1}^{k}  u_{ki} z_l  \nonumber \\
 \implies      & w_i - \sum_{l \in \{k+2,...,p\}} \bar{w}_{li} - \sum_{l =1}^{k} b_{ki} z_l   \geq w_i - \sum_{l \in \{k+1,...,p\}} \bar{w}_{li} -\sum_{l =1}^{k}  u_{ki} z_l. \label{eq:redundant_also_2} 
\end{align}
 \end{subequations}
 
Again, constraints (\ref{eq:redundant_2}) and (\ref{eq:redundant_also_2}) are redundant and can be eliminated, since $u_{ki} \geq b_{ki}$ for all $k \in [p]$. Through some minor rearranging, the resulting polyhedron is
\begin{subequations}
\begin{align}
Q^{proj}_{k+1} = Q^{const} \cap   \{  \bm{w}, y , \bar{\bm{w}}, \bar{\bm{y}}, \bm{z} ~~|~ & \sum_{l =1}^{k+1} b_{ki} z_l \leq w_i - \sum_{l \in \{k+2,...,p\}} \bar{w}_{li}  \\
  & \sum_{l =1}^{k+1}  u_{ki} z_l \geq w_i - \sum_{l \in \{k+2,...,p\}} \bar{w}_{li}   \\
&  u_{li} z_l \geq \bar{w}_{li}  \qquad~\forall l \in \{k+2,...,p\}\\
&  b_{li} z_l \leq \bar{w}_{li}  \qquad ~ \forall l \in \{k+2,...,p\} \}.
\end{align}
\end{subequations}

This proves the inductive step. After eliminating $ \bar{w}_{pi}$ from $Q^{proj}_{p}$, it should be clear that this results in
\begin{subequations}
\begin{align}
Q^{const} \cap   \{  \bm{w}, y , \bar{\bm{w}}, \bar{\bm{y}}, \bm{z} ~~|~ & \sum_{l =1}^{p} b_{ki} z_l \leq w_i   \\
& \sum_{l =1}^{p}  u_{ki} z_l \geq w_i  \}.
\end{align}
\end{subequations}

We can repeat the inductive procedure for the other features in the same manner. Finally $\bar{y}_l$ is eliminated for $l \in [p]$ by simple substitution of (\ref{union_score_app}) into (\ref{union_form_y_app}), and we arrive at formulation (\ref{projected_formulation}). The proof follows since the ideal property is preserved by projection.  \endproof

\section{Proof Lemma \ref{facet_def}}
\label{proof_facet_defining}
\proof{}
We show that (\ref{lower_bound_facet}) is facet defining. It can be proved that (\ref{upper_bound_facet}) is facet defining using the same argument. The dimension of this polyhedron is $p+d-1$. To show constraint (\ref{lower_bound_facet}) is facet defining, we need to find $p+d-1$ affinely independent points that satisfy $\sum_{l=1}^{p-1} b_{pi}+(b_{li}-b_{pi})z_l = w_i$.

Constraint (\ref{lower_bound_facet}) places bounds on dimension $i$ of $\bm{w}$. Without loss of generality, consider leaf $1$ and $i=d$. Define $\bm{\hat{w}}$ as a point on the interior of the leaf $\mathcal{L}_1$ with respect to dimensions $1,...,d-1$, but at the lower bound for dimension $d$, so that $\hat{w}_d=b_{1d}$. Define the point $\bm{q}^0=(\bm{\hat{w}},\bm{e}^1)$. This point satisfies (\ref{lower_bound_facet}) with equality.

Consider $d-1$ points $\bm{q}^i=(\bm{\hat{w}}+\epsilon \bm{e}^i,\bm{e}^1)$ for $i \in \{1,...,d-1\}$, where $\epsilon>0$ is chosen to be sufficiently small that $\bm{q}^i \in \mathcal{L}_1 $. A sufficiently small $\epsilon$ exists due to the fact that $\bm{\hat{w}}$ is on the interior with respect to dimensions $1,...,d-1$. These points still satisfy (\ref{lower_bound_facet}) with equality since $q^i_d=b_{1d}$. 

Consider $p-1$ points $\bm{\tilde{q}}^i=( \bm{\tilde{w}}^i,\bm{e}^i )$ for $i \in \{2,...,p-1\}$ and $\bm{\tilde{q}}^p=(\bm{\tilde{w}}^p,0)$, where $\bm{\tilde{w}}^i \in \mathcal{L}_i$ and $\tilde{w}^i_i=b_{id}$. Such a point exists because $\mathcal{L}_i$ is non-empty. 

We now need to show that these points are affinely independent. This can be proven by showing the following matrix is full rank:


\begin{align}
\begin{pmatrix} 
\bm{\tilde{q}}^2 - \bm{q}^0 \\
	\vdots \\ 
\bm{\tilde{q}}^{p} - \bm{q}^0 \\ 
\bm{q}^1 - \bm{q}^0 \\
	\vdots \\ 
\bm{q}^{d-1} - \bm{q}^0 \\ 
\end{pmatrix} = \begin{pmatrix} 
\bm{\tilde{w}}^2 - \bm{\hat{w}} & \bm{e}^2 - \bm{e}^1  \\
	 \vdots &  \vdots   \\ 
\bm{\tilde{w}}^p - \bm{\hat{w}} & \bm{0} - \bm{e}^1  \\
\epsilon\bm{e}^1 & \bm{0} \\
 \vdots  &	\vdots \\ 
\epsilon \bm{e}^{d-1} &  \bm{0} \\
\end{pmatrix}   .
\end{align}

If we shift the $p-2$ last columns to be the first $p-2$ columns and the $p-1$ to last column to $p-1$ from first, we end up with an upper diagonal matrix with nonzero entries on the diagonal, resulting in a matrix with full row rank. Since we applied only elementary operations to the original matrix, this also has full row rank. 
\begin{align}
\begin{pmatrix} 
\bm{\tilde{q}}^2 - \bm{q}^0 \\
\\
	\vdots \\ 
	\\
	\\
\bm{\tilde{q}}^{p} - \bm{q}^0 \\ 
\bm{q}^1 - \bm{q}^0 \\
	\vdots \\ 
	\\
\bm{q}^{d-1} - \bm{q}^0 \\ 
\end{pmatrix} = \begin{pmatrix} 
1 &  &  &  & -1 & \tilde{w}_1^2 - \hat{w}_1 &  & \dots & \tilde{w}_{d-1}^2 - \hat{w}_{d-1}  \\
	 & 1 & & & -1 & & &  & \\ 
	 & & \ddots & &\vdots & \vdots & & \ddots &  \\
	&  & & 1 & -1 &  & &  & \\
	&  & & & -1 & \tilde{w}_1^p - \hat{w}_1 &  & \dots  & \tilde{w}_{d-1}^2 - \hat{w}_{d-1}   \\
	&  & & & & \epsilon &  &   &  \\
	&  & & & &  &   &   &  \\
	&  & & & &  &  & \ddots &  \\
 &  & & & &  &  &  &  \epsilon  \\
\end{pmatrix}  .
\end{align}

This proves that the points were affinely independent and (\ref{lower_bound_facet}) is facet defining. 
\endproof

\section{Proof of Proposition \ref{above_below_prop}}
\label{proof_above_below_prop}
 \proof{}

Let $\bm{z},\bm{x}$ be any feasible solution to $ Q^{misic} \cap (\{0,1\}^p  \times \mathbb{R}^{1+p})$, where $\bm{x}$ is restricted to a binary lattice. We will show that $\bm{z},\bm{x}$ is feasible for $ Q^{expset} \cap (\{0,1\}^p  \times \mathbb{R}^{1+p}) $. 

For a given split $s$, suppose $x_{V(s)C(s)}=0$. Then  $x_{V(s')C(s')}=0$  for all $s'$ that have a lower threshold on the same variable, $C(s') \leq C(s), V(s') = V(s)$. This is due to constraint (\ref{bigger_than_cosntraint}), $x_{ij} \leq x_{ij+1}$, which enforces that $\bm{x}$ is a vector of 0's, followed by 1's. Therefore, combined with constraint (\ref{left_constraint}), all leaf variables $z_l$ are set to 0 for all leaves with thresholds less than $s$:
\begin{equation}
   \sum_{l \in \textbf{left}(s)} z_{l} \leq 0 ~~\forall s' ~s.t. ~ C(s') \leq C(s), V(s') = V(s). \label{all_left_0_equation}
\end{equation}

We analyze constraint (\ref{left_tight_constraint}) from $ Q^{expset}$ to check if it is satisfied. We begin by expanding the constraint out:
\begin{align*} \sum_{l \in \textbf{below}(s)} z_{l}  &= \sum_{l \in \textbf{left}(s)} z_{l} +  \sum_{l \in \textbf{below}(s_{-1}) \setminus \textbf{left}(s)} z_{l} \\
&= \sum_{l \in \textbf{left}(s)} z_{l} +  \sum_{l \in \textbf{left}(s_{-1}) \setminus \textbf{left}(s)} z_{l} 
+  \sum_{l \in \textbf{left}(s_{-2}) \setminus (\textbf{left}(s)\cup \textbf{left}(s_{-1}) )} z_{l} + ...  
\end{align*}
where $s_{-1}$ is the next threshold below $s$, such that $C(s)= C(s_{-1})+1$, and $s_{-2}$ is the next below that. This follows from the definition of the set $ \textbf{below}(s_{ij+1})= \textbf{below}(s_{ij}) \cup \textbf{left}(s_{ij+1})$. From equation (\ref{all_left_0_equation}),  all leaves with thresholds less than $s$ are set to 0, so:
\[ \sum_{l \in \textbf{below}(s)} z_{l} = 0 = x_{V(s)C(s)}. \]

Therefore, constraint (\ref{left_tight_constraint}) is satisfied. Furthermore, in this case (\ref{right_tight_constraint}) is trivially satisfied, since 
\[ \sum_{l \in \textbf{above}(s)} z_{l} = 1- x_{V(s)C(s)} = 1. \]

For the case where $x_{V(s)C(s)}=1$, the argument is very similar. In particular, since   $x_{V(s')C(s')}=1$, for all thresholds higher than $s$, it follows that
$$\sum_{l \in \textbf{right}(s)} z_{l} \leq 0 ~~\forall s' ~s.t. ~ C(s') \geq C(s), V(s') = V(s). $$

Analyzing constraint (\ref{right_tight_constraint}):
\begin{align*}
\sum_{l \in \textbf{above}(s)} z_{l}  =& \sum_{l \in \textbf{right}(s)} z_{l} +  \sum_{l \in \textbf{above}(s_{+1}) \setminus \textbf{right}(s)} z_{l} \\ 
= & \sum_{l \in \textbf{right}(s)} z_{l} +  \sum_{l \in \textbf{right}(s_{+1}) \setminus \textbf{right}(s)} z_{l} 
+  \sum_{l \in \textbf{right}(s_{+2}) \setminus (\textbf{right}(s)\cup \textbf{right}(s_{+1}) )} z_{l} ...  
\end{align*}
where $s_{+1},s_{+2}... $ are thresholds immediately above $s$. It follows that \[ \sum_{l \in \textbf{above}(s)} z_{l} = 0 = 1- x_{V(s)C(s)}.  \]

Again, (\ref{left_tight_constraint}) is trivially satisfied. We also have $Q^{expset} \subseteq Q^{misic}$ since:
\begin{align*}
 & \sum_{l \in \textbf{below}(s)} z_{l} \leq  x_{V(s)C(s)} \implies \sum_{l \in \textbf{left}(s)} z_{l} \leq  x_{V(s)C(s)} \\  
  & \sum_{l \in \textbf{above}(s)} z_{l} \leq  1- x_{V(s)C(s)} \implies \sum_{l \in \textbf{right}(s)} z_{l} \leq  1- x_{V(s)C(s)}.
\end{align*}

This occurs because $\textbf{left}  \subseteq \textbf{below}$ and $\textbf{right}  \subseteq \textbf{above}$.

 \endproof{}

\section{Proof of Proposition \ref{proof_tighter_subset}}
\label{ap_proof_tighter_subset}
 \proof{}

 For convenience, we recall the definition polyhedra $\tilde{Q}^{misic}(s,s')$ and $\tilde{Q}^{expset}(s,s'):$
  \begin{align*}
&&\tilde{Q}^{misic}(s,s')= \{ \bm{x},\bm{z} ~ |&~ \sum_{l \in \textbf{left}(s)} z_{l} \leq x_{V(s)C(s)}, ~  \sum_{l \in \textbf{right}(s)} z_{l} \leq 1- x_{V(s)C(s)}, \\
&&&~\sum_{l \in \textbf{left}(s')} z_{l} \leq x_{V(s')C(s')}, ~  \sum_{l \in \textbf{right}(s')} z_{l} \leq 1- x_{V(s')C(s')},\\
&&& x_{V(s)C(s)} \leq  x_{V(s')C(s')}\}\\
&&\tilde{Q}^{expset}(s,s')= \{ \bm{x},\bm{z} ~ |&~ \sum_{l \in \textbf{below}(s)} z_{l} \leq x_{V(s)C(s)}, ~  \sum_{l \in \textbf{above}(s)} z_{l} \leq 1- x_{V(s)C(s)}, \\
&&&~\sum_{l \in \textbf{below}(s')} z_{l} \leq x_{V(s')C(s')}, ~  \sum_{l \in \textbf{above}(s')} z_{l} \leq 1- x_{V(s')C(s')},\\
&&& x_{V(s)C(s)} \leq  x_{V(s')C(s')}\}.
\end{align*}
 
There are three cases that need to be examined: where split $s$ is a child of split $s'$, where $s'$ is a child of split $s$, and where neither is a child of the other because they are on different branches of the tree. Recall that in all cases, $s'$ is the split with a larger threshold. 

\begin{enumerate}
\item We start with the case where $s$ is a (left) child of split $s'$. An example of this occurs in Figure \ref{fig:biggs_not_ideal}. Take the solution $\bm{x}^{(1)},\bm{z}^{(1)}$ such that $\sum_{l \in \textbf{left}(s)} z_{l}^{(1)}=0, ~\sum_{l \in \textbf{left}(s)} z_{l}^{(1)}=0.5, ~\sum_{l \in  \textbf{right}(s)} z_{l}^{(1)}=0.5, ~\sum_{l \in  \textbf{right}(s')} z_{l}^{(1)}=0.5, ~x_{V(s)C(s)}^{(1)}=0.5, ~x_{V(s')C(s')}^{(1)}=0.5$.  By inspection, $\bm{x}^{(1)},\bm{z}^{(1)} \in \tilde{Q}^{misic}(s,s')$. This doesn't necessarily violate $\sum_{1=1}^p z_l^{(1)}=1$ since, $\textbf{left}(s') \supseteq \textbf{left}(s) \cup \textbf{right}(s)$.

Since $s'$ is the greater split, we have that $\textbf{above}(s) \supseteq \textbf{right}(s) \cup \textbf{right}(s')$. Furthermore, $\textbf{right}(s) \cap \textbf{right}(s')=\emptyset$, since $s$ is the left child of $s'$. It follows that the solution $\bm{x}^{(1)},\bm{z}^{(1)} \notin \tilde{Q}^{expset}(s,s')$ since 

$$ \sum_{l \in \textbf{above}(s)} z_{l}^{(1)}  \geq  \sum_{l \in \textbf{right}(s)} z_{l}^{(1)} + \sum_{l \in \textbf{right}(s')} z_{l}^{(1)} = 1 \implies \sum_{l \in \textbf{above}(s)} z_{l}^{(1)} \not\leq  1- x_{V(s)C(s)}^{(1)}. $$

\item We next examine the case where $s'$ is a (right) child of split $s$, which is very similar but included for completeness. Take the solution $\bm{x}^{(2)},\bm{z}^{(2)}$ such that $\sum_{l \in \textbf{left}(s)} z_{l}^{(2)}=0.5, ~\sum_{l \in \textbf{left}(s)} z_{l}^{(2)}=0.5, ~\sum_{l \in  \textbf{right}(s)} z_{l}^{(2)}=0.5, ~\sum_{l \in  \textbf{right}(s')} z_{l}^{(2)}=0, ~x_{V(s)C(s)}^{(2)}=0.5, ~x_{V(s')C(s')}^{(2)}=0.5$. By inspection, $\bm{x}^{(2)},\bm{z}^{(2)} \in \tilde{Q}^{misic}(s,s')$.

Since $s'$ is the greater split, we have that $\textbf{below}(s') \supseteq \textbf{left}(s) \cup \textbf{left}(s')$. Furthermore, $\textbf{left}(s) \cap \textbf{left}(s')=\emptyset$, since $s'$ is the right child of $s$. It follows that the solution $\bm{x}^{(2)},\bm{z}^{(2)} \notin \tilde{Q}^{expset}(s,s')$ since 

$$ \sum_{l \in \textbf{below}(s')} z_{l}^{(2)}  \geq  \sum_{l \in \textbf{left}(s)} z_{l}^{(2)} + \sum_{l \in \textbf{left}(s')} z_{l}^{(2)} = 1 \implies \sum_{l \in \textbf{below}(s')} z_{l}^{(2)} \not\leq  x_{V(s')C(s')}^{(2)}.$$

\item Finally, we examine the case where neither split is a child of the other, which is also very similar to the case above. An example of this occurs in Figure \ref{above_below_tree_example}. Take the solution $\bm{x}^{(3)},\bm{z}^{(3)}$ such that $\sum_{l \in \textbf{left}(s)} z_{l}^{(3)}=0.5, ~\sum_{l \in \textbf{left}(s)} z_{l}^{(3)}=0.5, ~\sum_{l \in  \textbf{right}(s)} z_{l}^{(3)}=0, ~\sum_{l \in  \textbf{right}(s')} z_{l}^{(3)}=0, ~x_{V(s)C(s)}^{(3)}=0.5, ~x_{V(s')C(s')}^{(3)}=0.5$. By inspection, $\bm{x}^{(3)},\bm{z}^{(3)} \in \tilde{Q}^{misic}(s,s')$.

Since $s'$ is the greater split, we have that $\textbf{below}(s') \supseteq \textbf{left}(s) \cup \textbf{left}(s')$. Furthermore, $\textbf{left}(s) \cap \textbf{left}(s')=\emptyset$, since neither node is a child of the other. It follows that the solution $\bm{x}^{(3)},\bm{z}^{(3)} \notin \tilde{Q}^{expset}(s,s')$ since 

$$ \sum_{l \in \textbf{below}(s')} z_{l}^{(3)}  \geq  \sum_{l \in \textbf{left}(s)} z_{l}^{(3)} + \sum_{l \in \textbf{left}(s')} z_{l}^{(3)} = 1 \implies \sum_{l \in \textbf{below}(s')} z_{l}^{(3)} \not\leq  x_{V(s')C(s')}^{(3)}.$$

For this case, there is also another fractional solution $\bm{x}^{(4)},\bm{z}^{(4)}$ such that $\sum_{l \in \textbf{left}(s)} z_{l}^{(4)}=0, ~\sum_{l \in \textbf{left}(s)} z_{l}^{(4)}=0, ~\sum_{l \in  \textbf{right}(s)} z_{l}^{(4)}=0.5, ~\sum_{l \in  \textbf{right}(s')} z_{l}^{(4)}=0.5, ~x_{V(s)C(s)}^{(4)}=0.5, ~x_{V(s')C(s')}^{(4)}=0.5$, which can be proven with a very similar argument. 
\end{enumerate}

\section{Proof of Proposition \ref{cuts_prop}}
\label{proof_cuts_prop}
 \proof{}
 
 Let $\bm{z},\bm{x}$ be any feasible solution to $ Q^{misic}\cap (\{0,1\}^p  \times \mathbb{R}^{1+p})$, where $\bm{x}$ is restricted to a binary lattice. We will show that $\bm{z},\bm{x}$ is feasible for $ Q^{elbow}\cap (\{0,1\}^p  \times \mathbb{R}^{1+p})$. 
 
 In particular, we will show that the $\bm{z},\bm{x}$ satisfies (\ref{left_tight_constraint_bend}). Consider two splits $s$ and $s'$ covered by the constraint (\ref{left_tight_constraint_bend}), where $s' \in \textbf{right\_parent}(s)$, that is $s'$ is above and to the right of $s$ in the tree. We will investigate the different feasible values for $x_{V(s)C(s)}, x_{V(s')C(s')}$, specifically $x_{V(s)C(s)}, x_{V(s')C(s')} \in \{(0,0),(0,1),(1,1)\}$.  Note that $x_{V(s)C(s)}=1, x_{V(s')C(s')}=0$ is not a feasible solution since it violates the constraint (\ref{bigger_than_cosntraint}), $x_{ij} \leq x_{ij+1}$.
 
 Suppose $x_{V(s)C(s)}=0 $ and $ x_{V(s')C(s')}=0$. From constraint (\ref{left_constraint}), $\sum_{l \in \textbf{left}(s')} z_{l} \leq 0$. However, since $s' \in \textbf{right\_parent}(S)$, then $\textbf{right}(s) \subset \textbf{left}(s')$. Therefore, $\sum_{l \in \textbf{right}(s)} z_{l} \leq 0$. As a result, (\ref{left_tight_constraint_bend}) is satisfied since:
\[  \sum_{l \in \textbf{right}(s)} z_{l}  \leq 0 =  x_{V(s')C(s')} - x_{V(s)C(s)}. \]
 
  Suppose $x_{V(s)C(s)}=0, x_{V(s')C(s')}=1$, then $  \sum_{l \in \textbf{right}(s)} z_{l}  \leq  x_{V(s')C(s')} - x_{V(s)C(s)} =1 $ is immediately satisfied.    Suppose $x_{V(s)C(s)}=1, x_{V(s')C(s')}=1$, then from (\ref{right_constraint}), $ \sum_{l \in \textbf{right}(s)} z_{l} \leq 1 - x_{V(s)C(s)} = 0  $,  therefore (\ref{left_tight_constraint_bend}) is also satisfied.

  It folows that $Q^{elbow} \subseteq Q^{misic}$, since $Q^{elbow}$ only has constraints added in addition to the constraints in $Q^{misic}$.

  \endproof{}
\section{Proof Theorem \ref{ideal_one_dimension}}
\label{proof_ideal_one_dimension}
 \proof{}
 
 For the case when $d=1$, constraint (\ref{right_tight_constraint}) can be rearranged:
 
\begin{align*}   \sum_{l \in \textbf{above}(s)} z_{l} \leq 1- x_{V(s)C(s)}  &\iff 1- \sum_{l \in \textbf{above}(s)} z_{l} \geq x_{V(s)C(s)} \\
&\iff  \sum_{l \in \textbf{below}(s)} z_{l} \geq x_{V(s)C(s)}.
\end{align*} 
 
 This is because in the single feature case, the sets $\textbf{above}(s)$ and $\textbf{below}(s)$ are complementary. Note that this isn't the case with multiple features, as generally there will be leaves in the tree that do not split on a feature.  This implies $\sum_{l \in \textbf{below}(s)} z_{l} = x_{V(s)C(s)} ~ \forall s ~\in  ~ \textbf{splits}(t)$ for $d=1$.
 
We now define the matrix $A$ by ordering the constraints in a specific way. We order the rows corresponding to variables $x_{j}$ from smallest to largest according to the size of the threshold to which they correspond. Furthermore, suppose leaves $z_{lt}$ are labelled in increasing order, so that $z_{1t}$ is the leaf corresponding to smallest threshold for tree $t$, while $z_{12t}$ is the next smallest. In this case, $A$ will take the following form:

\begin{equation}
\setstackgap{L}{1.0\baselineskip}
\fixTABwidth{T}
\parenMatrixstack{
1 &  &  &  &  &  &     &- 1&   &   &   &   &    \\
  &  &  & 1&  &  &     &   & -1&   &   &   &    \\
  &  &  & 1& 1&  &\hdots&   &   & -1&   &   &    \\
  &  &  & 1& 1& 1&     &   &   &   & -1&   &    \\
1 & 1&  &  &  &  &     &   &   &   &   &- 1&    \\
1 & 1& 1&  &  &  &     &   &   &   &   &   & -1   \\
 &  &  &\vdots&  &&\ddots   &   &   &\vdots  &   &   &   \\
  &  &  &  &  &  &     &  -1& 1&   &   &   &   \\
  &  &  &  &  &  &     &   &  -1& 1&   &   &   \\
  &  &  &  &  &  &\hdots&   &   &  -1& 1&   &   \\
  &  &  &  &  &  &     &   &   &   &  -1& 1&   \\
  &  &  &  &  &  &     &   &   &   &  &   -1& 1  }  .
\end{equation}


In the top right, we have the negative identity matrix, corresponding to each $x$. In the bottom right, we have the constraints $x_{ij} \leq x_{ij+1}$. In the top left, we have blocks of rows, each corresponding to leaves in a tree. In the example given, columns 1-3 correspond to leaves from tree 1 and columns 4-6 from tree 2. Due to the construction of the set $\textbf{below}$ where $ \textbf{below}(s_{ij+1})= \textbf{below}(s_{ij}) \cup \textbf{left}(s_{ij+1})$ for subsequent ordered splits in the tree, these rows have a lower triangular structure.

To prove any subset of this matrix is totally unimodular, we use the following lemma, originally from \cite{ghouila1962caracterisation}, presented as Theorem 19.3 in  \cite{schrijver1998theory}:

\begin{lem} (Ghouila-Houri)
Each collection of columns of A can be split into two parts so that the sum of the columns in one part minus the sum of the columns in the other part is a vector with entries only 0, + 1, and - 1
\end{lem}

To construct these sets, we allocate the first column (available in the subset of columns) of each tree to group 1. We then alternate through the remaining columns, assigning every second (available) column in the tree to group -1. This ensures that the sum of each row is 1 or 0 for the left columns corresponding to the $z_l$ leaf variables available in the subset, due to the consecutive ones property of the lower triangular matrix. We assign the remaining available columns (corresponding to $x$ variables) to group 1. The -1 from the identity matrix, if present in the subset, will reduce the sum to -1 or 0. For the lower half of the matrix corresponding to $x_{ij} \leq x_{ij+1}$, there is at most one 1, and one -1. Since these are assigned to the same group, the sum of the columns for these rows is either 0, + 1, or - 1. As a result, the total sum of all columns for all subsets is either 0, + 1, or - 1.  This assignment is illustrated below for a sample matrix where $\sigma$ corresponds to the group assignment.

\begin{equation*}
\setstackgap{L}{1.0\baselineskip}
\fixTABwidth{T}
\begin{blockarray}{cccccccccccccccc}
\sigma= &\blu 1 &\gre -1 &\blu 1 &\blu 1  & \gre -1 & \blu 1 &  \hdots  &\blu 1&\blu 1  &\blu 1  &\blu 1  &\blu 1  &\blu  1 & & A \sigma \bla \\
\begin{block}{c(ccccccccccccc)cc}
& &  &  &  &  &  &     & &   &   &   &   &   && \\
&\blu 1 &  &  &  &  &  &     &\blu  - 1&   &   &   &   &   && 0\\
&  &  &  &\blu 1&  &  &     &   &\blu  -1&   &   &   &   && 0  \\
&  &  &  &\blu 1&\gre 1&  &\hdots&   &   &\blu  -1&   &   &  &&  -1  \\
&  &  &  &\blu 1&\gre 1&\blu  1&     &   &   &   &\blu  -1&   &   &&  0 \\
&\blu 1 & \gre 1&  &  &  &  &     &   &   &   &   &\blu - 1&   &&  -1 \\
&\blu 1 & \gre 1& \blu1&  &  &  &     &   &   &   &   &   &\blu  -1&& 0   \\
A=& &  &  &\vdots&  &&\ddots   &   &   &\vdots  &&   &   &=& \vdots &   \\
&  &  &  &  &  &  &     &\blu   -1&\blu  1&   &   &   &  && 0\\
&  &  &  &  &  &  &     &   &\blu   -1&\blu  1&   &   &  && 0\\
&  &  &  &  &  &  &\hdots&   &   &\blu   -1& \blu 1&   & && 0\\
&  &  &  &  &  &  &     &   &   &   &\blu   -1&\blu  1&  && 0\\
&  &  &  &  &  &  &     &   &   &   &  &   \blu -1&\blu  1&&0  \\ 
\end{block}
\end{blockarray} 
\end{equation*}

 \endproof{}

\section{Proof of Lemma \ref{theorem_cuts_above_below}}
\label{sec:proof_theorem_cuts_above_below}

 \proof{}  

 We will prove the first statement, while the proof for the second statement is almost identical. To reduce the notation, assume the sets below are intercepted with $Q^{misic}$.
 \begin{align*}
 \sum_{l \in \textbf{below}(s')} & z_{l} \leq x_{V(s')C(s')}  \bigcap \sum_{l \in \textbf{above}(s)} z_{l}  \leq 1- x_{V(s)C(s)}   \\ &\implies
 \sum_{l \in \textbf{below}(s')} z_{l} + \sum_{l \in \textbf{above}(s)} z_{l}  \leq x_{V(s')C(s')}+ 1- x_{V(s)C(s)}  \\
 & \implies  \sum_{l \in \textbf{below}(s')} z_{l} + \sum_{l \in \textbf{above}(s)} z_{l} -1  \leq x_{V(s')C(s')} - x_{V(s)C(s)} \\
  & \implies  \sum_{l \in \textbf{below}(s')\setminus\textbf{right}(s)} z_{l} + \sum_{l \in \textbf{above}(s)\setminus \textbf{right}(s) } z_{l} + 2 \sum_{l \in \textbf{right}(s)} z_{l} -1 \\
  &\qquad \qquad \qquad \qquad \qquad \qquad \qquad  \leq x_{V(s')C(s')} - x_{V(s)C(s)} \\
    & \implies  1+  \sum_{l \in \textbf{right}(s)} z_{l} -1  \leq x_{V(s')C(s')} - x_{V(s)C(s)} \\
    & \implies   \sum_{l \in \textbf{right}(s)} z_{l}   \leq x_{V(s')C(s')} - x_{V(s)C(s)} .
 \end{align*}

The second-to-last implication follows because  $\textbf{right}(s) \subset \textbf{below}(s')$ and $\textbf{right}(s) \subset \textbf{above}(s)$. The last implication occurs because $\textup{\textbf{below}}(s') \cup \textup{\textbf{above}}(s) = p$, when combined with $\sum_{l=1}^p z_{l}=1$, we have that $ \sum_{l \in \textbf{below}(s')\setminus\textbf{right}(s)} z_{l} + \sum_{l \in \textbf{above}(s)\setminus \textbf{right}(s) } z_{l} + \sum_{l \in \textbf{right}(s)} z_{l} =1$.   
 
 \endproof{}

 \section{Problem size larger forests}
\label{sec:problem_size}

\begin{table*}[!htbp]
\centering
\caption{Problem sizes for instance with 5 features, depth 20}
\begin{tabular}{lllll}
\hline
\# trees            & method & constraints & binary variables & nonzeros \\
\hline
\multirow{4}{*}{1}  & \texttt{projected}   & 11          & 2766             & 27709    \\
                    & \texttt{misic}  & 8276        & 2761             & 54917    \\
                    & \texttt{bigM}   & 16560       & 8287             & 41398    \\
                    & \texttt{elbow}  & 8865        & 2761             & 61927    \\
                    \hline
\multirow{4}{*}{2}  & \texttt{projected}   & 22          & 5627             & 56857    \\
                    & \texttt{misic}   & 16873       & 5622            & 112953   \\
                    & \texttt{bigM}   & 33720       & 16869            & 84296    \\
                    & \texttt{elbow}  & 18038       & 5622            & 128593   \\
                    \hline
\multirow{4}{*}{4}  & \texttt{projected}   & 44          & 11312            & 114060   \\
                    & \texttt{misic}   & 34003       & 11307            & 225842   \\
                    & \texttt{bigM}   & 67818       & 33922            & 169537   \\
                    & \texttt{elbow}  & 36404       & 11307            & 254902   \\
                    \hline
\multirow{4}{*}{8}  & \texttt{projected}   & 88          & 22832            & 227507   \\
                    & \texttt{misic}   & 68964       & 22827            & 453909   \\
                    & \texttt{bigM}   & 136914      & 68478            & 342269   \\
                    & \texttt{elbow}  & 73692       & 22827            & 520911   \\
                    \hline
\multirow{4}{*}{16} & \texttt{projected}   & 176         & 45206            & 455015   \\
                    & \texttt{misic}   & 137322      & 45201            & 911007   \\
                    & \texttt{bigM}   & 271110      & 135592           & 677743   \\
                    & \texttt{elbow}  & 146789      & 45201            & 1032816  \\
                    \hline
\multirow{4}{*}{32} & \texttt{projected}   & 352         & 91990            & 924083   \\
                    & \texttt{misic}   & 282939      & 91985           & 1847111  \\
                    & \texttt{bigM}   & 551718      & 275928           & 1379231  \\
                    & \texttt{elbow}  & 302335      & 91985           & 2097640 
\end{tabular}
\label{table:problem_size}
\end{table*}

\section{Number of nodes in the branch and bound tree}

\blu 
To further analyze the computational behavior of different formulations, we examined the number of nodes explored in the branch-and-bound tree for the \texttt{concrete} and \texttt{winequalityred} datasets in the simulations presented in Section \ref{sec:solve_time_real_data}. The results show that formulations incorporating ordered binary feature variables, such as \texttt{elbow} and \texttt{misic}, tend to explore significantly fewer nodes compared to formulations like \texttt{bigM} and \texttt{projected}, which do not impose such an ordering. This difference is likely due to the branching decisions: in formulations with ordered binary feature variables, branching effectively prunes large portions of the search space by ``turning off" many other variables, thereby reducing the total number of nodes that need to be explored. However, while ordered binary representations may improve search efficiency, they can limit the flexibility of incorporating further constraints and may be less tight relative to the \texttt{projected} formulation. 

\begin{figure}
  \centering
  \begin{subfigure}{0.42\textwidth}
    \centering
	\includegraphics[width=\textwidth]{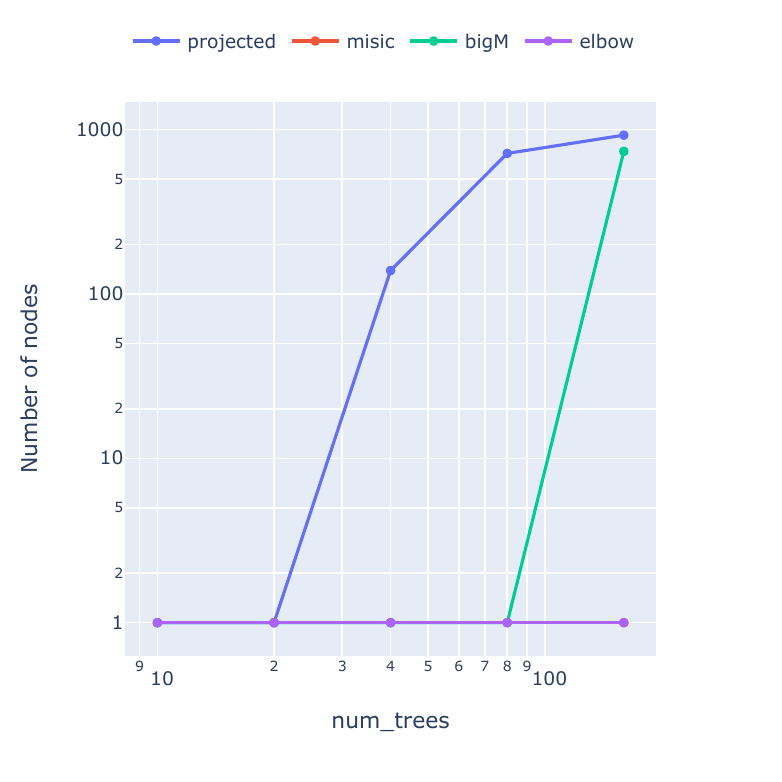}
	\caption{1 feature}
 \label{fig:concrete_1_feat_nodes}
  \end{subfigure}
    \begin{subfigure}{0.42\textwidth}
    \centering
	\includegraphics[width=\textwidth]{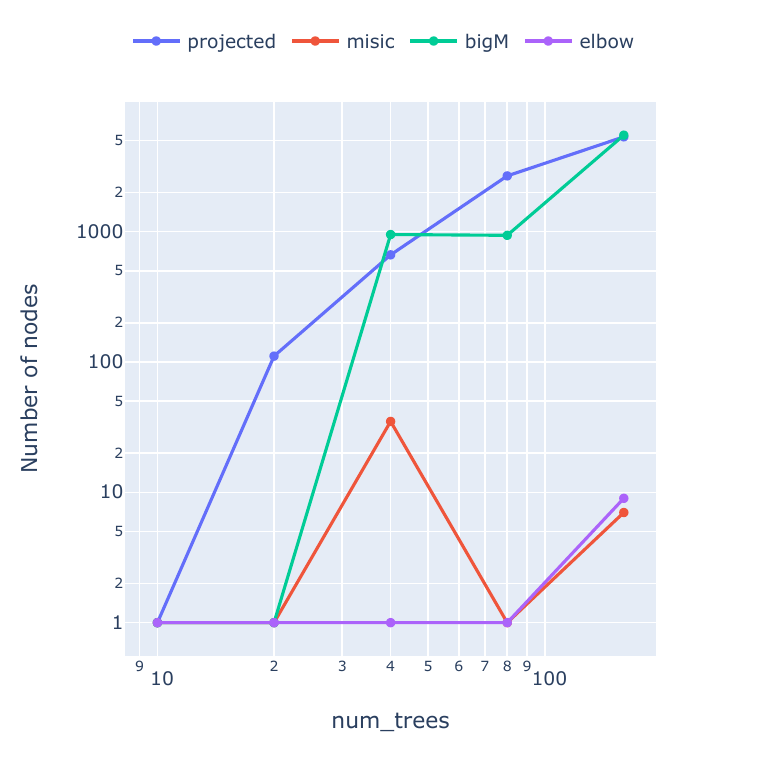}
	\caption{3 features}
  \end{subfigure}
    \begin{subfigure}{0.42\textwidth}
    \centering
	\includegraphics[width=\textwidth]{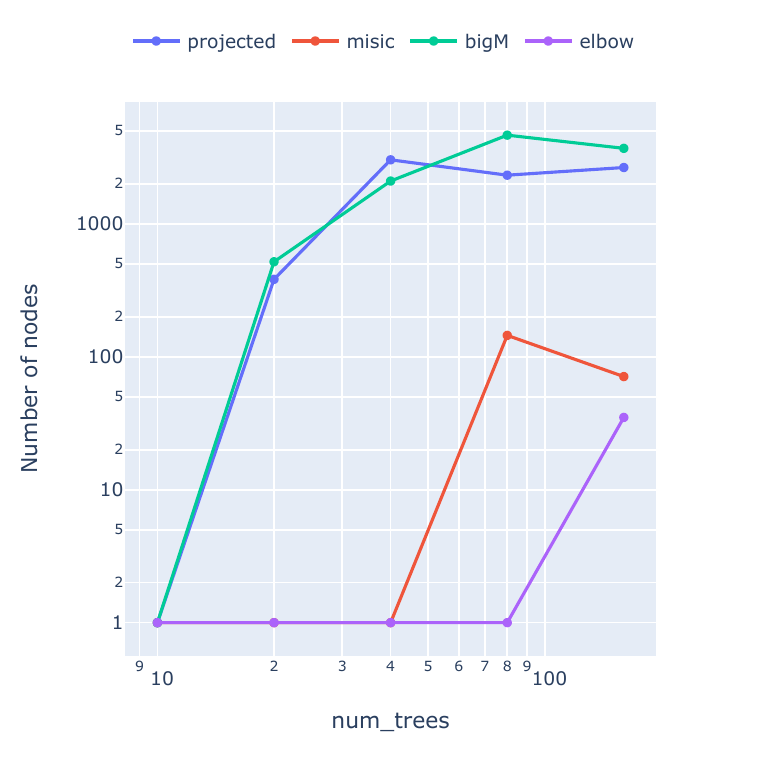}
	\caption{5 features}
  \label{fig:concrete_5_feat_nodes}
  \end{subfigure}
     \begin{subfigure}{0.42\textwidth}
    \centering
	\includegraphics[width=\textwidth]{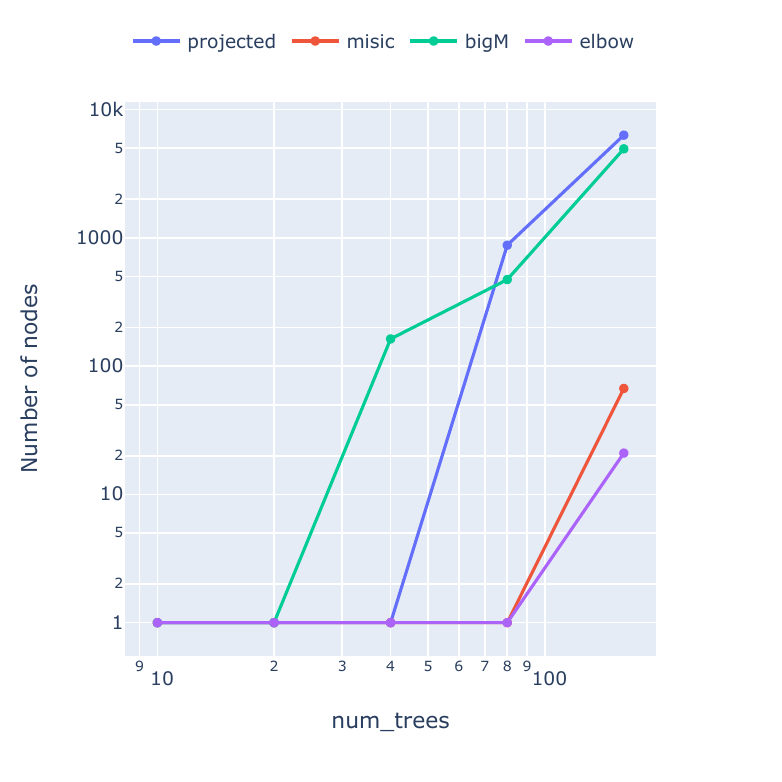}
	\caption{7 features}
  \end{subfigure}
  \caption{Number of nodes in the branch and bound tree for random forests of varying sizes \texttt{concrete} data}
  \label{fig:num_nodes_rf_concrete}
\end{figure}

\begin{figure}
  \centering
  \begin{subfigure}{0.42\textwidth}
    \centering
	\includegraphics[width=\textwidth]{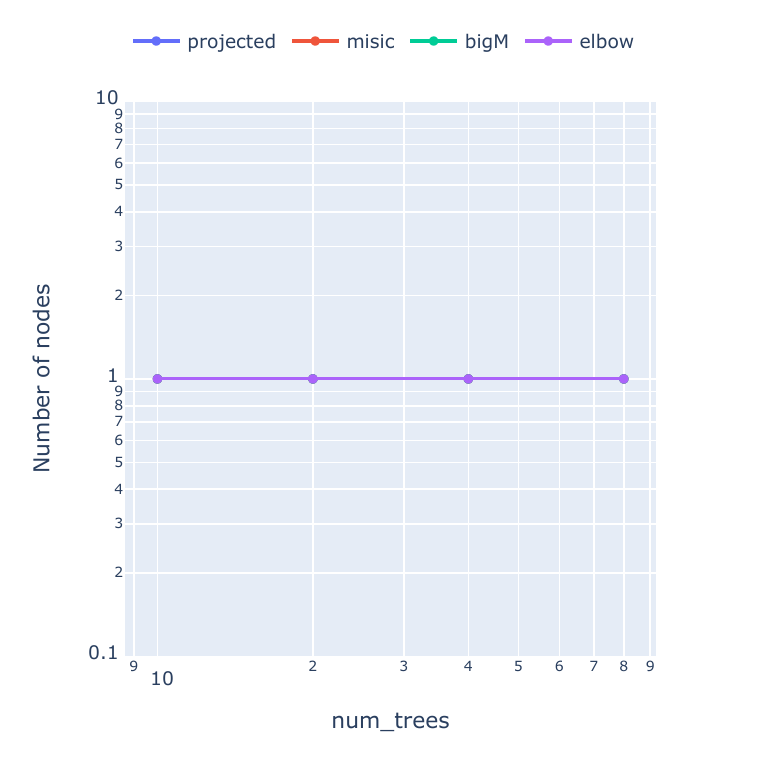}
	\caption{1 feature}
 \label{fig:wine_1_feat_nodes}
  \end{subfigure}
    \begin{subfigure}{0.42\textwidth}
    \centering
	\includegraphics[width=\textwidth]{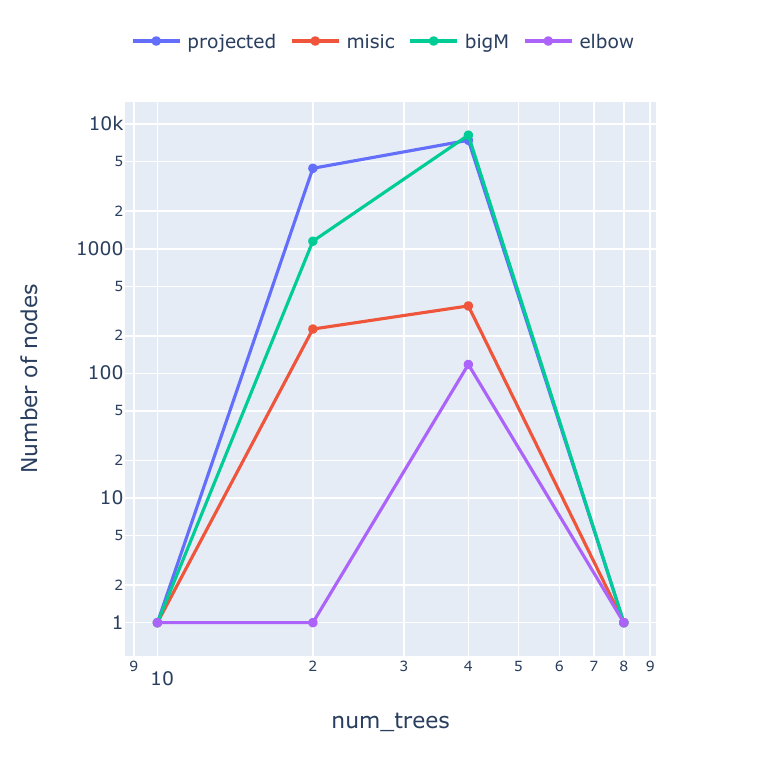}
	\caption{5 features}
  \end{subfigure}
    \begin{subfigure}{0.42\textwidth}
    \centering
	\includegraphics[width=\textwidth]{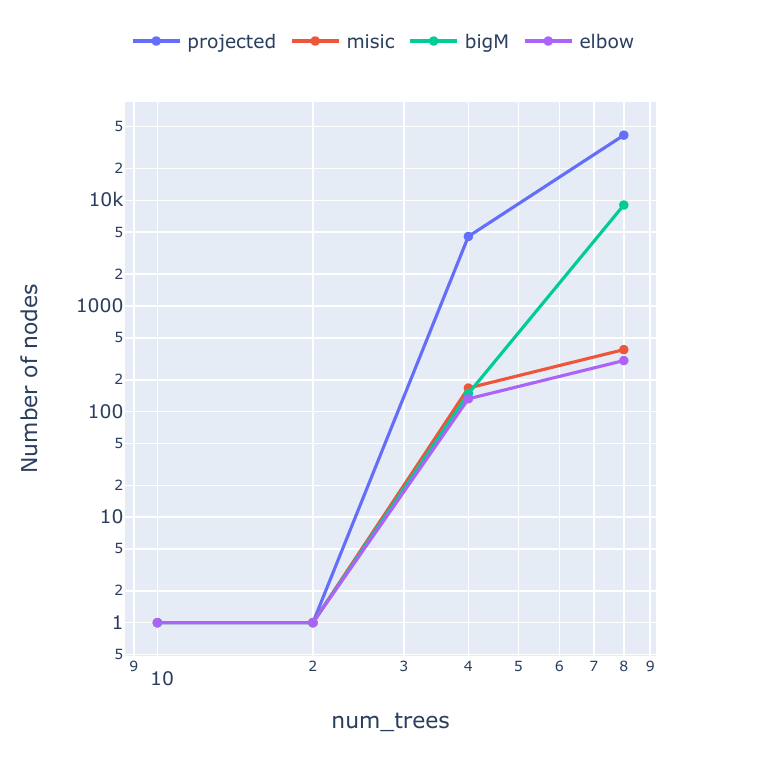}
	\caption{10 features}
  \label{fig:wine_5_feat_nodes}
  \end{subfigure}
  \caption{Number of nodes in the branch and bound tree for random forests of varying sizes \texttt{winequalityred} data}
  \label{fig:num_nodes_rf_wine}
\end{figure}

\bla